\pdfoutput=1
\documentclass{article}

\usepackage[paper=a4paper,left=30mm,right=30mm,top=25mm,bottom=25mm]{geometry}
\usepackage{amsmath,amssymb,amsthm,mathtools}
\usepackage{enumitem}
\usepackage{graphicx}
\usepackage{placeins}
\usepackage[numbers]{natbib}
\usepackage[colorlinks=true,linkcolor=blue,citecolor=blue,urlcolor=blue]{hyperref}

\newtheorem{theorem}{Theorem}[section]
\newtheorem{lemma}[theorem]{Lemma}
\newtheorem{proposition}[theorem]{Proposition}
\newtheorem{corollary}[theorem]{Corollary}
\theoremstyle{remark}
\newtheorem{remark}[theorem]{Remark}

\newcommand{\CC}{\mathcal{C}}
\newcommand{\Gam}{\Gamma}
\newcommand{\E}{\mathbb{E}}
\newcommand{\Prob}{\mathbb{P}}
\newcommand{\R}{\mathbb{R}}
\newcommand{\de}{\,\mathrm{d}}
\newcommand{\1}{\mathbf{1}}
\newcommand{\Leb}{\mathrm{Leb}}
\newcommand{\clip}[1]{\bigl[#1\bigr]_0^1}
\newcommand{\keywords}[1]{\par\medskip\noindent\textbf{Keywords:} #1}

\title{Exact regions for the symmetrized Blest coefficient}
\author{Marcus Rockel}
\date{\today}

\begin{document}

\maketitle
\begin{center}
  \small\textit{
    Department of Quantitative Finance,\\
    Institute for Economics, University of Freiburg,\\
    Rempartstr.~16, 79098 Freiburg, Germany,\\
    \texttt{marcus.rockel@finance.uni-freiburg.de}\\[2mm]
  }
\end{center}

\begin{abstract}
Spearman's rho treats all ranks alike, yet in many comparisons of two rankings an agreement at the top matters more than one at the bottom.
Blest's rank correlation builds this preference into Spearman's rho, and its symmetrized version $\eta$, introduced by Genest and Plante, does so without favoring either of the two rankings.
How far can such a top-weighted coefficient stray from other rank correlations?
We answer this question for Spearman's rho, Blomqvist's beta and Chatterjee's rank correlation $\xi$ by determining the exact set of values that $\eta$ and each of these coefficients can attain jointly.
Given Spearman's rho, $\eta$ can differ from it by at most $27/128$, and the extreme values of $\eta$ are attained only when the two rankings agree perfectly on the top ranks and are completely reversed on the bottom ranks, or vice versa.
Given Blomqvist's beta, $\eta$ obeys exactly the same sharp bounds as Spearman's rho.
For Chatterjee's rank correlation, both boundary curves are explicit and are attained by copulas with densities.
Here symmetrization has a surprising effect.
Compared with Blest's original coefficient, $\eta$ has a narrower range at fixed Spearman's rho, but near independence a wider range at fixed Chatterjee's rank correlation, where $|\eta|$ can come close to $\sqrt{11\xi/5}$ while Blest's coefficient stays below $\sqrt{32\xi/15}$.
The proofs combine explicit Kantorovich potentials from optimal transport with a convex relaxation that turns out to be exact.
\end{abstract}


\section{Introduction}

When two judges rank the same wines, or two search engines rank the same documents, a disagreement about the first places usually matters more than a disagreement about the last ones.
Spearman's rho, one of the most widely used measures of agreement between two rankings, is blind to this distinction, since it penalizes a swap of two neighboring ranks in the same way at the top and at the bottom of the list.
Blest \cite{blest2000rank} proposed a variant of Spearman's rho that weights the leading ranks of one ranking more heavily, at the price that the coefficient changes when the two rankings exchange their roles.
Genest and Plante \cite{genest2003blest} therefore introduced the average of Blest's coefficient over both orderings of the rankings.
It keeps the emphasis on the leading ranks, now of both rankings, and treats the two rankings alike.
As we point out below, this coefficient is also the population version of the weighted rank correlation of Pinto da Costa and Soares \cite{pintodacosta2005weighted}, which was designed with the same aim.

A top-weighted coefficient is useful only if it tells us something that Spearman's rho does not.
How far apart can the two coefficients be, and which patterns of agreement drive them apart?
More generally, how much information does the symmetrized coefficient share with other rank coefficients?
We answer these questions at the population level by computing exact regions, that is, the sets of all pairs of values that two coefficients attain jointly over all dependence structures.
The boundary of an exact region contains all sharp inequalities between the two coefficients, and the dependence structures on the boundary show what drives the coefficients furthest apart.

Let $\CC$ denote the class of bivariate copulas, i.e., of distribution functions on $[0,1]^2$ with uniform marginals, see \cite{nelsen2006introduction,durante2016principles}, and let $M(u,v)\coloneqq\min\{u,v\}$, $W(u,v)\coloneqq\max\{u+v-1,0\}$ and $\Pi(u,v)\coloneqq uv$ denote the comonotonicity, countermonotonicity and independence copulas.
For $C\in\CC$, the population version of Blest's coefficient is $\nu(C)=24\int_0^1\!\int_0^1(1-u)\,C(u,v)\de u\de v-2$, see \cite[Sec.~3]{genest2003blest}.
As in \cite{blest2000rank}, small values of $u$ correspond to the first places of a ranking, so that the weight $1-u$ emphasizes the leading ranks.
With the transpose $C^\top(u,v)\coloneqq C(v,u)$, the symmetrized Blest coefficient is
\begin{equation}
  \label{eq:eta-def}
  \eta(C)\coloneqq\frac{\nu(C)+\nu(C^\top)}{2}
  =12\int_0^1\!\int_0^1(2-u-v)\,C(u,v)\de u\de v-2,
\end{equation}
see \cite[Sec.~4]{genest2003blest}.
Both coefficients equal $1$ at $M$, $-1$ at $W$ and $0$ at $\Pi$.
The coefficient $\eta$ is invariant under transposition, but, as observed in \cite[Sec.~4]{genest2003blest}, it violates the reflection axiom of a measure of concordance in the sense of Scarsini \cite{scarsini1984measures}.

The role of the weight in~\eqref{eq:eta-def} becomes transparent when $\eta$ is compared with Spearman's rho $\rho(C)=12\int_0^1\!\int_0^1C(u,v)\de u\de v-3$, see \cite[Sec.~5.1.2]{nelsen2006introduction}.
For $(U,V)\sim C$, both coefficients are affine functions of an average squared rank difference,
\[
  \rho(C)=1-6\,\E\bigl[(U-V)^2\bigr],
  \qquad
  \eta(C)=1-6\,\E\bigl[(U-V)^2(2-U-V)\bigr],
\]
see Section~\ref{sec:preliminaries}.
Thus $\eta$ is Spearman's rho with the squared rank differences weighted by $2-U-V$.
The weight has mean $1$, decreases from $2$ at the leading ranks $U=V=0$ to $0$ at the trailing ranks $U=V=1$, and is symmetric in the two rankings.
For two rankings of $n$ objects, in which the $i$th object has the ranks $R_i$ and $Q_i$, the corresponding weights $(n-R_i+1)+(n-Q_i+1)$ define the weighted rank correlation $r_W$ of \cite{pintodacosta2005weighted}, and $\eta$ is the population version of $r_W$.

We compare $\eta$ with Spearman's rho, with Blomqvist's beta $\beta(C)=4C(\frac12,\frac12)-1$, introduced in \cite{blomqvist1950measure}, and with Chatterjee's rank correlation
\begin{equation}
  \label{eq:xi-def}
  \xi(C)\coloneqq6\int_0^1\!\int_0^1\bigl(\partial_1C(u,v)\bigr)^2\de u\de v-2.
\end{equation}
The latter is the limit of the rank statistic of Chatterjee \cite{chatterjee2020} and appeared earlier as a copula-based measure of regression dependence in \cite{dette2013copula}.
For $(U,V)\sim C$, it vanishes if and only if $U$ and $V$ are independent, and it equals one if and only if $V$ is almost surely a measurable function of $U$, see \cite[Thm.~1.1]{chatterjee2020}.
In contrast to $\rho$, $\beta$ and $\eta$, it measures the strength of a directed dependence rather than its sign.

The exact region of Kendall's tau and Spearman's rho was determined by Schreyer, Paulin and Trutschnig \cite{schreyer2017exact}, who showed that the classical inequalities between the two coefficients describe the region only partially.
Since then, many of the pairwise regions among Spearman's rho, Kendall's tau, Spearman's footrule, Gini's gamma and Blomqvist's beta have been determined, see, e.g., \cite{kokolbukovsek2019relation,kokolbukovsek2023exact,kokolbukovsek2024exact}.
For Chatterjee's rank correlation, the exact region with Spearman's rho was found in \cite{ansari2026exact}.
For Blest's coefficient, Tschimpke \cite{tschimpke2025exact} determined the region with Spearman's rho within the class of bivariate extreme-value copulas, the region with Chatterjee's rank correlation was established in \cite{rockel2026exact}, and the regions of $\nu$ with $\rho$ and with $\eta$ over all copulas were computed in \cite{rockel2026blest}.
To the best of our knowledge, the $(\nu,\eta)$-region of \cite{rockel2026blest} is the only exact region involving the symmetrized coefficient $\eta$ so far.

Our first result answers the opening question for Spearman's rho.
For $s\in[0,1]$, let $S_s$ denote the straight shuffle of min, cf.~\cite{mikusinski1992shuffles}, that is comonotone on $[0,1-s]$ and countermonotone on $[1-s,1]$, i.e., the copula of $(U,T_s(U))$ for a uniform random variable $U$ and $T_s(t)\coloneqq t$ for $t\le1-s$ and $T_s(t)\coloneqq2-s-t$ for $t>1-s$.
Its survival copula $\widehat S_s$ is countermonotone on $[0,s]$ and comonotone on $[s,1]$, see Figure~\ref{fig:shuffles} in Section~\ref{sec:beta}.
With
\begin{equation}
  \label{eq:Psi-def}
  \Psi(r)\coloneqq1-2^{-1/3}(1-r)^{4/3},\qquad r\in[-1,1],
\end{equation}
the region reads as follows.

\begin{theorem}[Exact $(\rho,\eta)$-region]
\label{thm:rho-eta}
It holds that
\[
  \bigl\{(\rho(C),\eta(C)):C\in\CC\bigr\}
  =\bigl\{(r,e)\in[-1,1]^2:\ 2r-\Psi(r)\le e\le\Psi(r)\bigr\}.
\]
For $r\in(-1,1)$ and $s\coloneqq(\frac{1-r}{2})^{1/3}$, the upper bound $e=\Psi(r)$ is attained by $S_s$ and by no other copula, and the lower bound $e=2r-\Psi(r)$ is attained by $\widehat S_s$ and by no other copula.
\end{theorem}

\begin{figure}[t!]
  \centering
  \includegraphics[width=\textwidth]{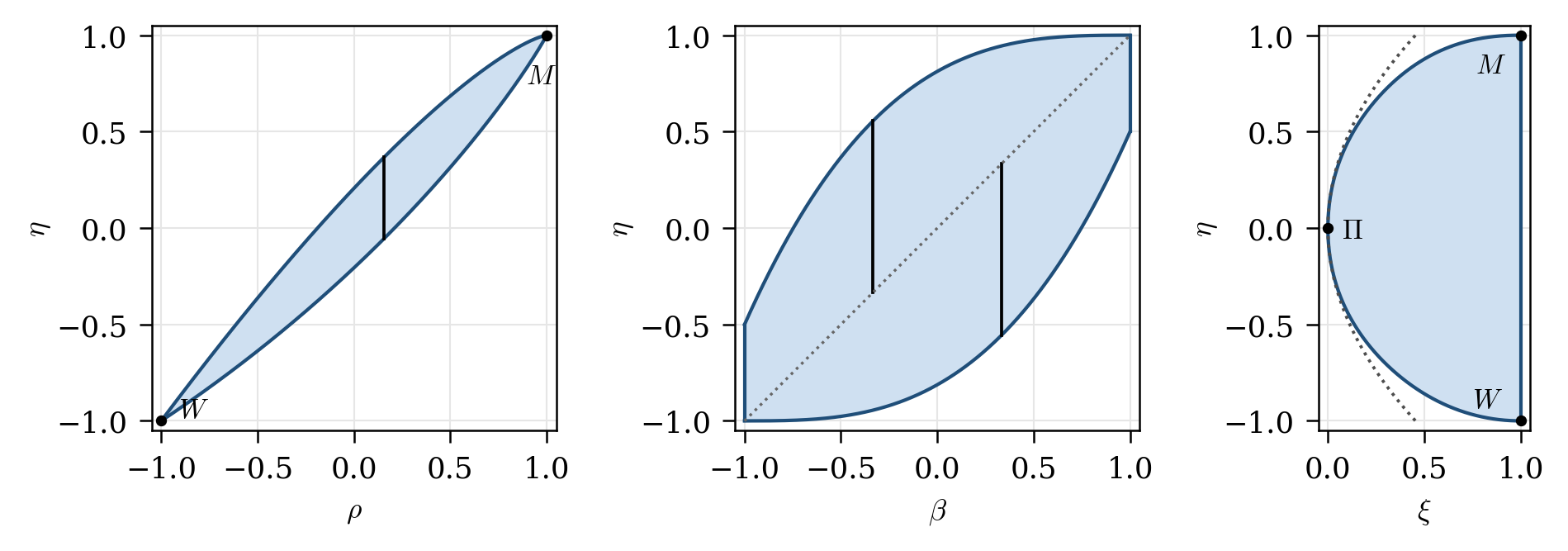}
  \caption{Exact regions of the symmetrized Blest coefficient $\eta$ with Spearman's rho (left, Theorem~\ref{thm:rho-eta}), Blomqvist's beta (middle, Theorem~\ref{thm:beta-eta}) and Chatterjee's rank correlation (right, Theorem~\ref{thm:xi-eta}), drawn from the closed-form boundaries.
  Left: the vertical segment marks the maximal deviation $|\eta-\rho|=\frac{27}{128}$ at $\rho=\frac{5}{32}$, and the comparison with Blest's coefficient is shown in Figure~\ref{fig:eta-nu}.
  Middle: the vertical segments mark the maximal deviation $|\eta-\beta|=\frac89$, attained at $\beta=\mp\frac13$.
  Right: the dotted curves $\pm\sqrt{11\xi/5}$ from Proposition~\ref{prop:sqrt} bound the region and are asymptotically attained as $\xi\to0$.}
  \label{fig:regions}
\end{figure}

The theorem is illustrated in Figure~\ref{fig:regions}, and its extremal copulas have a transparent meaning.
Under $S_s$, the two rankings agree perfectly on the leading ranks $[0,1-s]$ and are completely reversed on the trailing ranks $[1-s,1]$, so that all the disagreement compatible with $\rho=r$ is placed where $\eta$ is least sensitive to it.
Under $\widehat S_s$, the roles of the leading and the trailing ranks are exchanged.
For instance, if Spearman's rho equals $\frac12$, then $\eta$ can take every value between $2^{-5/3}\approx0.31$ and $1-2^{-5/3}\approx0.69$, and no other value.
In general, the fibers of the region are symmetric about the diagonal, and the theorem yields the sharp inequality $|\eta(C)-\rho(C)|\le\frac{27}{128}\approx0.21$, see Corollary~\ref{cor:eta-rho-deviation}.
Compared with the constant $\frac14$ for Blest's coefficient from \cite{rockel2026blest}, symmetrization thus brings the coefficient slightly closer to Spearman's rho.

Blomqvist's beta records only how likely the two variables are to lie on the same side of their medians.
It therefore constrains $\eta$ only weakly, and in fact exactly as weakly as it constrains Spearman's rho.
With
\begin{equation}
  \label{eq:rpm-def}
  r_-(b)\coloneqq-1+\tfrac{3}{16}(1+b)^3,
  \qquad
  r_+(b)\coloneqq1-\tfrac{3}{16}(1-b)^3,
  \qquad b\in[-1,1],
\end{equation}
the exact $(\beta,\rho)$-region is $\{(b,r):r_-(b)\le r\le r_+(b)\}$, see \cite{nelsen2001bounds,kokolbukovsek2019relation}, and our second result shows that the same holds for $\eta$ in place of $\rho$.

\begin{theorem}[Exact $(\beta,\eta)$-region]
\label{thm:beta-eta}
It holds that
\[
  \bigl\{(\beta(C),\eta(C)):C\in\CC\bigr\}
  =\bigl\{(b,e)\in[-1,1]^2:\ r_-(b)\le e\le r_+(b)\bigr\}.
\]
For every $b\in[-1,1]$, both bounds are attained by the symmetric shuffles of min from \eqref{eq:T-maps} below.
\end{theorem}

In particular, $|\eta(C)-\beta(C)|\le\frac89$ for every copula $C$, with equality at $\beta=\mp\frac13$, see Corollary~\ref{cor:eta-beta-deviation}.
The extremal copulas are again symmetric shuffles of min, and they are extremal for Spearman's rho as well.

Chatterjee's rank correlation behaves differently.
It measures how close $V$ is to being a function of $U$, and it is a quadratic rather than an affine functional of the copula.
Accordingly, the extremal dependence structures are no longer singular shuffles but absolutely continuous copulas.
Write $\clip{t}\coloneqq\min\{\max\{t,0\},1\}$ for $t\in\R$.
For $c>0$ and $v\in(0,1)$, there are unique numbers $a_c(v),b_c(v)\in\R$ such that
\begin{equation}
  \label{eq:ab-def}
  \int_0^1\clip{\frac{(2-u-v)^2-a_c(v)}{c}}\de u
  =v
  =\int_0^1\clip{\frac{b_c(v)-(2-u-v)^2}{c}}\de u,
\end{equation}
and by Lemma~\ref{lem:families} below, there are unique copulas $C_c^+$ and $C_c^-$ with
\begin{equation}
  \label{eq:families}
  \begin{gathered}
  \partial_1C_c^+(u,v)=\clip{\frac{(2-u-v)^2-a_c(v)}{c}},\\
  \partial_1C_c^-(u,v)=\clip{\frac{b_c(v)-(2-u-v)^2}{c}}
  \end{gathered}
\end{equation}
for almost all $(u,v)\in[0,1]^2$.
At every level $v$, the conditional distribution functions of these copulas are thus clamped quadratic functions of the rank sum $u+v$, see Figure~\ref{fig:xi-extremals} in Section~\ref{sec:xi}.

As the parameter $c$ runs through $(0,\infty)$, the points $(\xi(C_c^\pm),\eta(C_c^\pm))$ trace the two boundary curves of the $(\xi,\eta)$-region.
These curves are explicit.
In Section~\ref{sec:xi}, we define functions $x_\pm,e_\pm:(0,\infty)\to\R$ in closed form, see~\eqref{eq:up-large}--\eqref{eq:lo-small}.
They are built from rational functions, logarithms and the arctangent, evaluated at $\sqrt c$ or at the root of an explicit quartic equation with coefficients depending on $\sqrt c$, and Proposition~\ref{prop:closed} shows that $x_\pm(c)=\xi(C_c^\pm)$ and $e_\pm(c)=\eta(C_c^\pm)$.

\begin{theorem}[Exact $(\xi,\eta)$-region]
\label{thm:xi-eta}
Let $x_\pm$ and $e_\pm$ be the functions given in~\eqref{eq:up-large}--\eqref{eq:lo-small}.
Then $x_+$ and $x_-$ are continuous and strictly decreasing bijections from $(0,\infty)$ onto $(0,1)$.
Define $\eta_\pm:[0,1]\to[-1,1]$ by $\eta_\pm(0)\coloneqq0$, $\eta_\pm(1)\coloneqq\pm1$ and
\[
  \eta_\pm\bigl(x_\pm(c)\bigr)\coloneqq e_\pm(c),\qquad c>0.
\]
Then
\[
  \bigl\{(\xi(C),\eta(C)):C\in\CC\bigr\}
  =\bigl\{(x,e)\in[0,1]\times\R:\ \eta_-(x)\le e\le\eta_+(x)\bigr\}.
\]
The function $\eta_+$ is continuous, strictly increasing and concave, and $\eta_-$ is continuous, strictly decreasing and convex.
For every $c>0$, the copula $C_c^+$ is the only copula $C$ with $\xi(C)\le x_+(c)$ and $\eta(C)=e_+(c)$, and $C_c^-$ is the only copula $C$ with $\xi(C)\le x_-(c)$ and $\eta(C)=e_-(c)$.
\end{theorem}

A comparison with Blest's coefficient reveals an unexpected effect of symmetrization.
Since $\rho(C^\top)=\rho(C)$, both $\nu(C)$ and $\nu(C^\top)$ lie in the fiber of the $(\rho,\nu)$-region at $\rho(C)$, and so does their average $\eta(C)$.
Symmetrization can therefore only shrink the region with Spearman's rho.
Indeed, the $(\rho,\eta)$-region is a proper subset of the $(\rho,\nu)$-region, and the two boundaries meet only at $(-1,-1)$ and $(1,1)$, see Remark~\ref{rem:rho-nu}.
One might expect the same for Chatterjee's rank correlation, all the more so as the $(\xi,\eta)$-region closely resembles the $(\xi,\nu)$-region determined in \cite{rockel2026exact}, whose upper boundary differs from that of the $(\xi,\eta)$-region by less than $2\cdot10^{-3}$.
This expectation turns out to be false.
Near independence, $\eta$ can deviate from zero further than $\nu$ at the same value of $\xi$.
Indeed, we show the sharp inequality
\[
  |\eta(C)|\le\sqrt{\tfrac{11}{5}\,\xi(C)},
  \qquad C\in\CC,
\]
which is strict unless $C=\Pi$ and asymptotically attained as $\xi(C)\to0$, see Proposition~\ref{prop:sqrt}.
The same argument gives the sharp constants $2$ and $\frac{32}{15}$ in place of $\frac{11}{5}$ for Spearman's rho and for Blest's coefficient.
This is possible because $\xi$, unlike $\rho$, is not invariant under transposition, so that $\nu(C^\top)$ need not lie in the fiber of the $(\xi,\nu)$-region at $\xi(C)$.
Moreover, the $(\xi,\nu)$-region is symmetric about $\nu=0$, while the $(\xi,\eta)$-region is not, see Remark~\ref{rem:comparison}.

The proofs rest on the reflected coordinates $X\coloneqq1-U$ and $Z\coloneqq1-V$ for $(U,V)\sim C$, in which the weights $1-u$ and $2-u-v$ of $\nu$ and $\eta$ become $x$ and $x+z$, and $\eta(C)=6\,\E[XZ(X+Z)]-2$, see \cite[Lem.~2.1]{rockel2026blest}.
For Theorem~\ref{thm:rho-eta}, the supporting-line problems of the convex $(\rho,\eta)$-region are optimal transport problems with costs of the form $xz(x+z+\mu)$, which we solve with explicit piecewise cubic Kantorovich potentials, cf.~\cite{villani2009optimal}.
Theorem~\ref{thm:beta-eta} follows from a layer-cake decomposition of $\E[f(X)Z]$ combined with the bathtub principle.
For Theorem~\ref{thm:xi-eta}, we write both coefficients in terms of the conditional distribution functions $z\mapsto\Prob(Z\le z\mid X=x)$ and drop the requirement that these functions be nondecreasing in the threshold $z$.
The relaxed problem then decouples into one convex problem for each threshold, and its solutions turn out to be nondecreasing in $z$ nevertheless.
The relaxation is therefore exact, and the solutions are the kernels of the copulas $C_c^\pm$.

The paper is organized as follows.
Section~\ref{sec:preliminaries} collects the representation of the coefficients in reflected coordinates and the shuffles used as extremizers.
Sections~\ref{sec:rho}, \ref{sec:beta} and~\ref{sec:xi} prove Theorems~\ref{thm:rho-eta}, \ref{thm:beta-eta} and~\ref{thm:xi-eta}, respectively, together with the sharp deviation bounds and the comparison with Blest's coefficient.
Section~\ref{sec:discussion} records the elementary regions of $\eta$ with the tail dependence coefficients and discusses the pairs that remain open.
The appendix contains the derivation of the closed-form boundary curves of Theorem~\ref{thm:xi-eta}.

\section{Notation and preliminaries}
\label{sec:preliminaries}

For $(U,V)\sim C\in\CC$, we work with the reflected coordinates $X\coloneqq1-U$ and $Z\coloneqq1-V$, which are again uniform on $[0,1]$.
In these coordinates, a large value of $X$ signals a leading rank of the first variable, and the weights $1-u$ and $2-u-v$ in the definitions of $\nu$ and $\eta$ turn into $x$ and $x+z$, see Lemma~\ref{lem:coupling} below.
Conversely, every coupling of two uniform random variables $(X,Z)$ arises in this way from exactly one copula, and we identify $\CC$ with the set $\Gam$ of these couplings whenever convenient.
For $\gamma\in\Gam$, we write $\E_\gamma$ for expectations under $\gamma$ and drop the subscript when the coupling is clear.

The set $\CC$ is convex and compact with respect to uniform convergence, see \cite[Thm.~1.4.5 and Thm.~1.7.7]{durante2016principles}.
Lebesgue measure on $[0,1]$ is denoted by $\Leb$.
The following representation is \cite[Lem.~2.1]{rockel2026blest}, and the formula for $\beta$ is immediate from $C(\frac12,\frac12)=\Prob(X\ge\frac12,Z\ge\frac12)$.

\begin{lemma}[Coupling representation, \cite{rockel2026blest}]
\label{lem:coupling}
If $(X,Z)$ follows the coupling $\gamma\in\Gam$ associated with $C\in\CC$, then
\begin{equation}
  \label{eq:coupling-moments}
  \begin{gathered}
  \rho(C)=12\,\E[XZ]-3,
  \qquad
  \beta(C)=4\,\gamma\bigl([\tfrac12,1]^2\bigr)-1,\\
  \nu(C)=12\,\E[X^2Z]-2,
  \qquad
  \nu(C^\top)=12\,\E[XZ^2]-2,\\
  \eta(C)=6\,\E[XZ(X+Z)]-2.
  \end{gathered}
\end{equation}
In particular, $\rho$, $\nu$, $\eta$ and $\beta$ are affine functions of the copula.
\end{lemma}

Expanding the squares and using $\E[X^2]=\E[Z^2]=\frac13$ and $\E[X^3]=\E[Z^3]=\frac14$ gives
\[
  \E\bigl[(X-Z)^2\bigr]=\tfrac23-2\,\E[XZ],
  \qquad
  \E\bigl[(X-Z)^2(X+Z)\bigr]=\tfrac12-\E\bigl[XZ(X+Z)\bigr].
\]
Since $X-Z=V-U$ and $X+Z=2-U-V$, Lemma~\ref{lem:coupling} therefore yields the representations $\rho(C)=1-6\,\E[(U-V)^2]$ and $\eta(C)=1-6\,\E[(U-V)^2(2-U-V)]$ stated in the introduction.

The survival copula $\widehat C(u,v)\coloneqq u+v-1+C(1-u,1-v)$ of $C$ is the copula of $(1-U,1-V)$, and the reflection $C^\perp(u,v)\coloneqq u-C(u,1-v)$ is the copula of $(U,1-V)$.
By \cite[Lem.~2.2]{rockel2026blest} and the invariance of $\beta$ under the survival transformation, see \cite{scarsini1984measures}, we have
\begin{equation}
  \label{eq:survival}
  \rho(\widehat C)=\rho(C),
  \qquad
  \beta(\widehat C)=\beta(C),
  \qquad
  \eta(\widehat C)=2\rho(C)-\eta(C),
  \qquad C\in\CC.
\end{equation}
The survival transformation thus reflects the fiber $\{\eta(C):\rho(C)=r\}$ about the value $r$.
As noted in \cite[Sec.~2]{rockel2026blest}, the representation \eqref{eq:coupling-moments} shows that the reflection acts on $\eta$ as
\begin{equation}
  \label{eq:reflection}
  \eta(C^\perp)=-\rho(C)+\frac{\nu(C^\top)-\nu(C)}{2}.
\end{equation}
In particular, $\eta(C^\perp)=-\rho(C)\neq-\eta(C)$ for every symmetric copula $C$ with $\eta(C)\ne\rho(C)$, so that no transformation of this kind negates $\eta$ together with $\rho$.

Accordingly, the $(\rho,\eta)$-region below is not symmetric about the origin, see Remark~\ref{rem:rho-nu}, and the fibers of the $(\xi,\eta)$-region are not symmetric about zero, see Proposition~\ref{prop:closed}.

All explicit extremizers in Sections~\ref{sec:rho} and~\ref{sec:beta} are shuffles of min, cf.~\cite{mikusinski1992shuffles}.
If a Borel map $T:[0,1]\to[0,1]$ preserves $\Leb$, then the law of $(X,T(X))$ with $X$ uniform is a coupling in $\Gam$.
Its copula is symmetric whenever $T(T(x))=x$ for almost every $x$.
The shuffles $S_s$ from the introduction correspond to the coupling supported on the graph of $z=s-x$ for $x\in[0,s]$ and $z=x$ for $x\in(s,1]$, and their survival copulas $\widehat S_s$ to the graph of $z=x$ for $x\in[0,1-s]$ and $z=2-s-x$ for $x\in(1-s,1]$.
Both maps are involutions up to null sets, so $S_s$ and $\widehat S_s$ are symmetric, and $S_0=\widehat S_0=M$ and $S_1=\widehat S_1=W$.

\begin{lemma}[Coefficients along the shuffle families]
\label{lem:shuffle-values}
For every $s\in[0,1]$,
\[
  \rho(S_s)=\rho(\widehat S_s)=1-2s^3,
  \qquad
  \eta(S_s)=1-2s^4,
  \qquad
  \eta(\widehat S_s)=1-4s^3+2s^4.
\]
In particular, $s\mapsto\rho(S_s)$ is a decreasing bijection from $[0,1]$ onto $[-1,1]$, and $\eta(S_{s(r)})=\Psi(r)$ for $s(r)\coloneqq(\frac{1-r}{2})^{1/3}$ with $\Psi$ from~\eqref{eq:Psi-def}.
\end{lemma}

\begin{proof}
Integrating along the two linear pieces of the support of the coupling of $S_s$ gives
\begin{align*}
  \E[XZ]&=\int_0^sx(s-x)\de x+\int_s^1x^2\de x=\frac{s^3}{6}+\frac{1-s^3}{3},\\
  \E[XZ(X+Z)]&=s\int_0^sx(s-x)\de x+2\int_s^1x^3\de x=\frac{s^4}{6}+\frac{1-s^4}{2},
\end{align*}
so that \eqref{eq:coupling-moments} yields $\rho(S_s)=1-2s^3$ and $\eta(S_s)=1-2s^4$.
The values for $\widehat S_s$ follow from~\eqref{eq:survival}, and eliminating $s$ via $2s^3=1-r$ gives $\eta(S_{s(r)})=1-2(\frac{1-r}{2})^{4/3}=\Psi(r)$.
\end{proof}

\section{\texorpdfstring{The exact $(\rho,\eta)$-region}{The exact (rho,eta)-region}}
\label{sec:rho}

We first explain the strategy.
By Lemma~\ref{lem:coupling}, the map $C\mapsto(\rho(C),\eta(C))$ is affine, so the $(\rho,\eta)$-region is the image of the convex set $\CC$ under an affine map and hence convex.
A convex set is determined by its supporting lines, and by the survival symmetry~\eqref{eq:survival}, the lower boundary of the region is the reflection of the upper boundary about the diagonal.
It therefore suffices to solve the supporting-line problems
\[
  \sup\bigl\{\eta(C)+\lambda\rho(C):C\in\CC\bigr\},\qquad\lambda\in\R.
\]
By Lemma~\ref{lem:coupling}, the objective can be written as
\begin{equation}
  \label{eq:support-functional}
  \eta(C)+\lambda\rho(C)=6\,\E\bigl[XZ(X+Z+2\lambda)\bigr]-2-3\lambda,
\end{equation}
so each supporting-line problem is an optimal transport problem: we maximize the expectation of the cost $xz(x+z+2\lambda)$ over all couplings of two uniform random variables.

The shape of this cost suggests the optimal couplings.
Its mixed derivative is
\[
  \partial_x\partial_z\bigl\{xz(x+z+2\lambda)\bigr\}=2(x+z)+2\lambda,
\]
which changes sign on the unit square when $\lambda\in(-2,0)$.
Where it is positive, the cost is supermodular and rewards comonotone pairings, and where it is negative, near the origin, it is submodular and rewards countermonotone pairings.
In the reflected coordinates, the origin corresponds to the trailing ranks of both variables.
We are thus led to couplings that are countermonotone on a block near the origin and comonotone elsewhere, which are exactly the couplings of the shuffles $S_s$.
We verify their optimality with Kantorovich duality, cf.~\cite{villani2009optimal}.
A coupling is optimal as soon as there is a function $\varphi$ with $\varphi(x)+\varphi(z)\ge c(x,z)$ everywhere and with equality on the support of the coupling.
Since the cost is symmetric, we may use the same potential $\varphi$ in both variables.

For $s\in(0,1]$, let $\lambda_s\coloneqq-\frac{4s}{3}$, so that the cost becomes
\[
  c_s(x,z)\coloneqq xz\Bigl(x+z-\frac{8s}{3}\Bigr),
\]
and recall that the coupling of $S_s$ is supported on
\[
  \Sigma_s\coloneqq\bigl\{(x,s-x):x\in[0,s]\bigr\}\cup\bigl\{(t,t):t\in[s,1]\bigr\}.
\]
Along $\Sigma_s$, a potential with equality must satisfy the envelope condition $\varphi_s'(x)=\partial_xc_s(x,z)$ at the point $(x,z)\in\Sigma_s$ above $x$.
For $x>s$, this gives $\varphi_s'(x)=3x^2-\frac{8s}{3}x$, and for $x<s$, it gives $\varphi_s'(x)=-x^2+\frac{8s}{3}x-\frac{5s^2}{3}$.
Integrating and fixing the constants by continuity at $s$ and by equality at the point $(0,s)\in\Sigma_s$ leads to
\[
  \varphi_s(t)\coloneqq
  \begin{cases}
    -\frac{1}{3}t^3+\frac{4s}{3}t^2-\frac{5s^2}{3}t+\frac{s^3}{3}, & 0\le t\le s,\\[1mm]
    t^3-\frac{4s}{3}t^2, & s<t\le1.
  \end{cases}
\]
Both branches take the value $-\frac{s^3}{3}$ at $t=s$, so $\varphi_s$ is continuous.
The choice $\lambda_s=-\frac{4s}{3}$ is the one for which the two branches fit together in this way.

\begin{proposition}[Kantorovich potentials]
\label{prop:dual}
Let $s\in(0,1]$.
For all $(x,z)\in[0,1]^2$,
\begin{equation}
  \label{eq:dual-inequality}
  \varphi_s(x)+\varphi_s(z)\ge c_s(x,z),
\end{equation}
with equality if and only if $(x,z)\in\Sigma_s$.
\end{proposition}

\begin{proof}
Write $\Delta(x,z)\coloneqq\varphi_s(x)+\varphi_s(z)-c_s(x,z)$.
Since $\Delta$ is symmetric in $x$ and $z$, it suffices to consider three cases, according to the position of $x$ and $z$ relative to $s$.
Each of the following factorizations is verified by expanding both sides.
\begin{enumerate}[label=(\roman*),leftmargin=2em]
\item \emph{$s\le x,z\le1$.}
Here both arguments lie on the upper branch of $\varphi_s$, and
\[
  3\Delta(x,z)=(x-z)^2\bigl(3(x+z)-4s\bigr).
\]
This is nonnegative because $x+z\ge2s>\frac{4s}{3}$, and it vanishes if and only if $x=z$, that is, on the diagonal part of $\Sigma_s$.
\item \emph{$0\le x,z\le s$.}
Here both arguments lie on the lower branch, and
\[
  3\Delta(x,z)=(2s-x-z)(x+z-s)^2.
\]
This is nonnegative because $x+z\le2s$.
It vanishes if and only if $x+z=s$, which is the antidiagonal part of $\Sigma_s$, or $x=z=s$, and the point $(s,s)$ belongs to $\Sigma_s$ as well.
\item \emph{$0\le x\le s\le z\le1$.}
Substituting $x=s-v$ and $z=s+w$ with $v\in[0,s]$ and $w\in[0,1-s]$ gives $3\Delta(x,z)=g_v(w)$ with
\[
  g_v(w)\coloneqq(v+w)^2(2s-3v+3w)+v(s-2v)^2.
\]
The derivative $g_v'(w)=(v+w)(4s-3v+9w)$ is nonnegative because $v\le s$, so $g_v$ is nondecreasing on $[0,1-s]$.
Moreover, $g_v(0)=v(s-v)^2\ge0$.
Hence $\Delta\ge0$, and equality forces $w=0$ and $v\in\{0,s\}$, i.e., $(x,z)\in\{(s,s),(0,s)\}\subseteq\Sigma_s$.
\qedhere
\end{enumerate}
\end{proof}

\begin{proof}[Proof of Theorem~\ref{thm:rho-eta}]
Fix $s\in(0,1]$ and a coupling $\gamma\in\Gam$.
Integrating~\eqref{eq:dual-inequality} with respect to $\gamma$ and using that both marginals of $\gamma$ are uniform, we obtain
\[
  \E_\gamma\bigl[c_s(X,Z)\bigr]\le\E_\gamma\bigl[\varphi_s(X)+\varphi_s(Z)\bigr]=2\int_0^1\varphi_s(t)\de t.
\]
Equality holds for the coupling of $S_s$, which is concentrated on the equality set $\Sigma_s$.
In view of~\eqref{eq:support-functional}, this means that
\begin{equation}
  \label{eq:support-line}
  \eta(C)+\lambda_s\rho(C)\le\eta(S_s)+\lambda_s\rho(S_s),
  \qquad C\in\CC,\ s\in(0,1].
\end{equation}
We now derive the four parts of the theorem from this inequality.

\emph{Upper boundary.}
Let $C\in\CC$ with $\rho(C)=r$.
For $r\in(-1,1)$, we have $s(r)\in(0,1)$, and $\rho(S_{s(r)})=r$ by Lemma~\ref{lem:shuffle-values}.
Inequality~\eqref{eq:support-line} with $s=s(r)$ therefore reduces to the bound $\eta(C)\le\eta(S_{s(r)})$, and $\eta(S_{s(r)})=\Psi(r)$.
For $r=1$, the bound $\eta(C)\le1=\Psi(1)$ holds for every copula.
For $r=-1$, we take $s=1$ in~\eqref{eq:support-line} and obtain $\eta(C)\le\eta(S_1)=-1=\Psi(-1)$.

\emph{Lower boundary.}
The survival copula $\widehat C$ also satisfies $\rho(\widehat C)=r$, so the upper boundary applies to it.
By~\eqref{eq:survival}, we obtain $\eta(C)=2r-\eta(\widehat C)\ge2r-\Psi(r)$.

\emph{Attainment and fibers.}
For $r\in[-1,1]$ and $s=s(r)$, the copulas $S_s$ and $\widehat S_s$ attain the two bounds by Lemma~\ref{lem:shuffle-values}.
Since $\rho$ and $\eta$ are affine, the mixtures $(1-t)S_s+t\widehat S_s$, $t\in[0,1]$, all have Spearman's rho $r$, and their values of $\eta$ fill the fiber $[2r-\Psi(r),\Psi(r)]$.

\emph{Uniqueness.}
Let $r\in(-1,1)$, and let $C$ satisfy $\rho(C)=r$ and $\eta(C)=\Psi(r)$.
Then equality holds in~\eqref{eq:support-line} for $s=s(r)$, and the associated coupling $\gamma$ satisfies
\[
  \int_{[0,1]^2}\bigl(\varphi_s(x)+\varphi_s(z)-c_s(x,z)\bigr)\de\gamma(x,z)=0.
\]
The integrand is nonnegative, so $\gamma$ is concentrated on its zero set, which is $\Sigma_s$ by Proposition~\ref{prop:dual}.
For $x\in[0,s)$, the only point of $\Sigma_s$ with first coordinate $x$ is $(x,s-x)$, and for $x\in(s,1]$, it is $(x,x)$.
Since the first marginal of $\gamma$ is uniform, $\gamma$ is the coupling of $S_s$, and hence $C=S_s$.
Uniqueness of the lower extremizer follows because the survival transformation is an involution on $\CC$.
\end{proof}

Along the fiber at $\rho=r$, the deviation $\eta-\rho$ ranges over the symmetric interval $[-\omega(r),\omega(r)]$, where
\[
  \omega(r)\coloneqq\Psi(r)-r
\]
is the half-width of the fiber.
Maximizing $\omega$ gives the sharp deviation between the two coefficients.

\begin{corollary}[Deviation from Spearman's rho]
\label{cor:eta-rho-deviation}
Every copula $C$ satisfies $|\eta(C)-\rho(C)|\le\frac{27}{128}$.
The two extremes $\pm\frac{27}{128}$ occur only at $\rho=\frac{5}{32}$, where they are attained precisely by $S_{3/4}$ and $\widehat S_{3/4}$, respectively.
Moreover, $|\eta(C)|\le1-2^{-1/3}\approx0.2063$ whenever $\rho(C)=0$, with equality precisely for $S_{2^{-1/3}}$ and $\widehat S_{2^{-1/3}}$.
\end{corollary}

\begin{proof}
The derivative $\omega'(r)=\frac43\,2^{-1/3}(1-r)^{1/3}-1$ is strictly decreasing, so $\omega$ is strictly concave.
Its unique maximizer is given by $(1-r)^{1/3}=\frac34\,2^{1/3}$, that is, by $r=\frac{5}{32}$, and then $s(r)=\frac34$.
The maximal value is
\[
  \omega\bigl(\tfrac{5}{32}\bigr)=\Psi\bigl(\tfrac{5}{32}\bigr)-\tfrac{5}{32}=\tfrac{47}{128}-\tfrac{20}{128}=\tfrac{27}{128}.
\]
The second claim is the fiber of Theorem~\ref{thm:rho-eta} at $r=0$, where $\Psi(0)=1-2^{-1/3}$ and $s(0)=2^{-1/3}$.
The attainment statements follow from the uniqueness part of Theorem~\ref{thm:rho-eta}.
\end{proof}

\begin{figure}[t!]
  \centering
  \includegraphics[width=\textwidth]{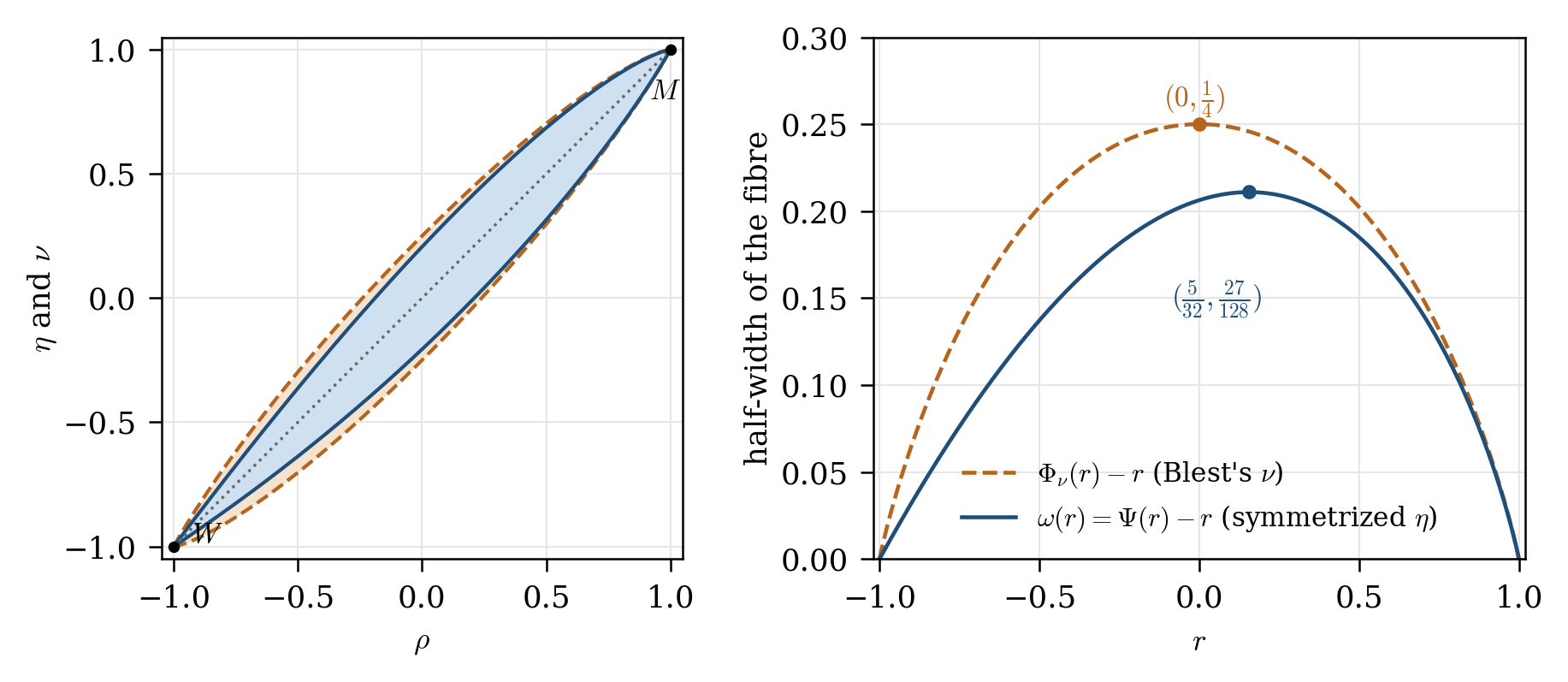}
  \caption{Comparison of the symmetrized coefficient $\eta$ with Blest's coefficient $\nu$ at fixed Spearman's rho.
  Left: the exact $(\rho,\eta)$-region of Theorem~\ref{thm:rho-eta} (blue) inside the exact $(\rho,\nu)$-region of \cite{rockel2026blest} (orange, dashed boundary), with the diagonal dotted.
  Right: the half-widths of the fibers at $\rho=r$, that is, the largest possible deviations $\omega(r)=\Psi(r)-r$ of $\eta$ and $\Phi_\nu(r)-r$ of $\nu$ from $\rho$.
  Symmetrization lowers the sharp deviation from $\frac14$, attained at $\rho=0$, to $\frac{27}{128}$, attained at $\rho=\frac{5}{32}$, and it makes the region asymmetric about the origin.
  All curves are drawn from their closed forms.}
  \label{fig:eta-nu}
\end{figure}

We close this section by comparing the region with that of Blest's coefficient itself, see Figure~\ref{fig:eta-nu}.

\begin{remark}[Comparison with Blest's coefficient]
\label{rem:rho-nu}
\leavevmode
\begin{enumerate}[label=(\alph*),leftmargin=2em]
\item
By \cite[Thm.~1.1]{rockel2026blest}, the fiber of the exact $(\rho,\nu)$-region at $\rho=r$ is the interval $[2r-\Phi_\nu(r),\Phi_\nu(r)]$, where
\[
  \Phi_\nu(r)=1-\tfrac34(1-r)^{4/3}\ \text{ for }r\ge0,
  \qquad
  \Phi_\nu(r)=1+2r-\tfrac34(1+r)^{4/3}\ \text{ for }r\le0.
\]
Since $\rho(C^\top)=\rho(C)$, both $\nu(C)$ and $\nu(C^\top)$ lie in this fiber, and so does their average $\eta(C)$.
The inclusion of the $(\rho,\eta)$-region in the $(\rho,\nu)$-region is therefore automatic.
It is strict except at $M$ and $W$, as an elementary comparison of $\Psi$ and $\Phi_\nu$ on $(-1,1)$ shows, which for $r\ge0$ reduces to $2^{-1/3}>\frac34$.
\item
The half-width of the $(\rho,\nu)$-region is maximal at $\rho=0$, with value $\frac14$, see \cite{rockel2026blest}.
For $\eta$, the maximum moves to $\rho=\frac{5}{32}$ and drops to $\frac{27}{128}\approx0.211$.
\item
The half-width $\omega$ is not an even function.
For instance, $\omega(\frac12)\approx0.185$ and $\omega(-\frac12)\approx0.137$, so that $\eta$ can deviate from $\rho$ more strongly under moderately positive than under moderately negative dependence.
The $(\rho,\nu)$-region, in contrast, is symmetric about the origin, and the asymmetry of the $(\rho,\eta)$-region reflects the failure of point symmetry in~\eqref{eq:reflection}.
\end{enumerate}
\end{remark}

\section{\texorpdfstring{The exact $(\beta,\eta)$-region}{The exact (beta,eta)-region}}
\label{sec:beta}

Blomqvist's beta only fixes the value of the copula at the point $(\frac12,\frac12)$.
In the reflected coordinates, prescribing $\beta(C)=b$ amounts to prescribing the mass $q\coloneqq\frac{1+b}{4}\in[0,\frac12]$ of the upper right quadrant
\[
  Q\coloneqq[\tfrac12,1]^2,
\]
because $\beta(C)=4\gamma(Q)-1$ by~\eqref{eq:coupling-moments}.
Since $\beta$ is invariant under transposition, our strategy is to bound $\nu(C)$ and $\nu(C^\top)$ separately at fixed $\beta$ and then to average.
The bounds for $\eta$ obtained in this way are sharp if the extremal copulas for $\nu$ are symmetric, because $\eta=\nu$ for symmetric copulas.
By~\eqref{eq:coupling-moments}, bounding $\nu$ means bounding $\E[X^2Z]$ under the constraint $\gamma(Q)=q$, a rearrangement problem with a constraint on one quadrant.

\begin{figure}[t!]
  \centering
  \includegraphics[width=\textwidth]{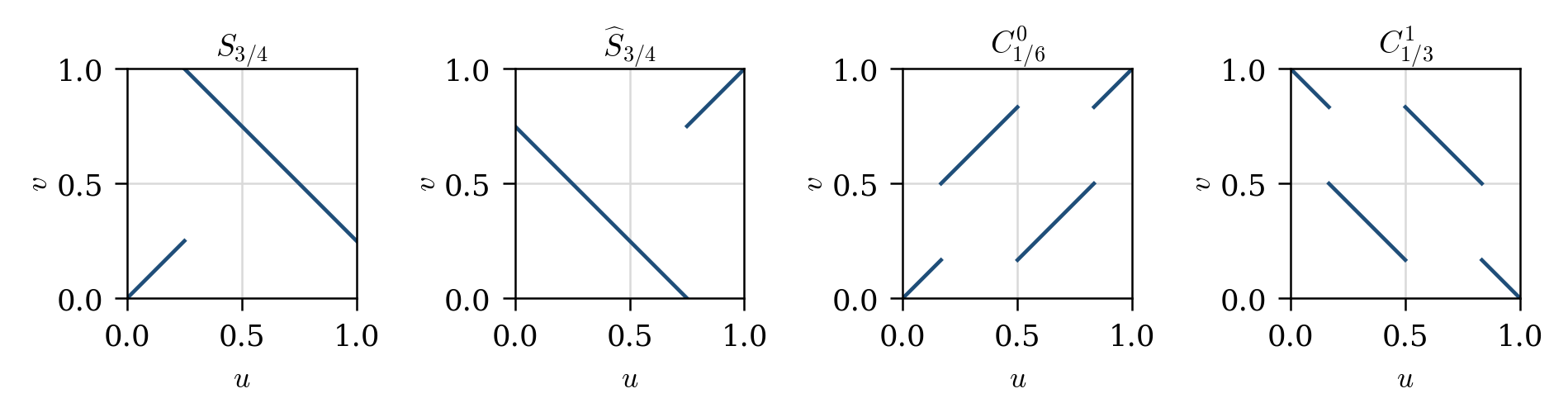}
  \caption{Supports of the extremal shuffles of Theorems~\ref{thm:rho-eta} and~\ref{thm:beta-eta} in the original coordinates $(u,v)$.
  The shuffle $S_{3/4}$ attains the maximal deviation $\eta-\rho=\frac{27}{128}$ and its survival copula $\widehat S_{3/4}$ the minimal deviation $-\frac{27}{128}$, both at $\rho=\frac{5}{32}$.
  The copulas $C^0_{1/6}$ and $C^1_{1/3}$ attain the extreme deviations $\eta-\beta=\frac89$ at $\beta=-\frac13$ and $\eta-\beta=-\frac89$ at $\beta=\frac13$, respectively.
  All four copulas are symmetric, so that $\eta=\nu$ for each of them.}
  \label{fig:shuffles}
\end{figure}

The extremal couplings are induced by the following two measure-preserving maps of $[0,1]$, both involutions up to null sets:
\begin{equation}
  \label{eq:T-maps}
  T^0_q(x)\coloneqq
  \begin{cases}
    x, & x\in[0,q]\cup[1-q,1],\\
    x+\frac12-q, & x\in(q,\frac12],\\
    x-\frac12+q, & x\in(\frac12,1-q),
  \end{cases}
  \qquad
  T^1_q(x)\coloneqq1-T^0_{1/2-q}(x).
\end{equation}
Explicitly, $T^1_q(x)=1-x$ for $x\in[0,\frac12-q]\cup[\frac12+q,1]$, $T^1_q(x)=1-q-x$ for $x\in(\frac12-q,\frac12]$, and $T^1_q(x)=1+q-x$ for $x\in(\frac12,\frac12+q)$.
We let $\gamma^0_q$ and $\gamma^1_q$ be the laws of $(X,T^0_q(X))$ and $(X,T^1_q(X))$ for a uniform random variable $X$, and $C^0_q$ and $C^1_q$ the associated copulas.
Both couplings satisfy $\gamma(Q)=q$, and both copulas are symmetric, see Figure~\ref{fig:shuffles}.
In words, $T^0_q$ is the identity on the outer blocks $[0,q]$ and $[1-q,1]$ and exchanges the two inner blocks by translation, so that the largest values of each half of $[0,1]$ are sent comonotonically into the upper half.
The map $T^1_q$ does the opposite: it sends the smallest values of each half countermonotonically into the upper half.

\begin{proposition}[Constrained rearrangement]
\label{prop:constrained}
Let $q\in[0,\frac12]$ and let $f:[0,1]\to\R$ be bounded and nondecreasing.
For every $\gamma\in\Gam$ with $\gamma(Q)=q$,
\[
  \E_{\gamma^1_q}[f(X)Z]\le\E_\gamma[f(X)Z]\le\E_{\gamma^0_q}[f(X)Z].
\]
\end{proposition}

\begin{proof}
We first prove the upper bound in four steps and then deduce the lower bound from it.

\emph{Step 1: layer-cake decomposition.}
Since $Z=\int_0^1\1\{Z>t\}\de t$, Fubini's theorem gives
\begin{equation}
  \label{eq:layer-cake}
  \E_\gamma[f(X)Z]=\int_0^1\E_\gamma\bigl[f(X)\1\{Z>t\}\bigr]\de t,
\end{equation}
and we bound each layer separately.
Fix $t\in[0,1)$ and let $m_t(x)$ be a version of the conditional probability $\gamma(Z>t\mid X=x)$ with values in $[0,1]$.
Then
\[
  \E_\gamma\bigl[f(X)\1\{Z>t\}\bigr]=\int_0^1f(x)m_t(x)\de x,
  \quad
  \int_0^1m_t(x)\de x=\gamma(Z>t)=1-t.
\]

\emph{Step 2: the constraint on each layer.}
Split $[0,1]$ into the halves $I_-\coloneqq[0,\frac12)$ and $I_+\coloneqq[\frac12,1]$, and let $\alpha\coloneqq\int_{I_+}m_t(x)\de x$ be the mass that the layer places above $I_+$, so that $\int_{I_-}m_t(x)\de x=1-t-\alpha$.
The quadrant constraint restricts $\alpha$.
If $t\ge\frac12$, then $\{Z>t\}\subseteq\{Z\ge\frac12\}$, and hence $\alpha\le\gamma(Q)=q$.
If $t<\frac12$, then $\{Z>t\}\supseteq\{Z\ge\frac12\}$, and hence $1-t-\alpha\ge\gamma(I_-\times I_+)=\frac12-q$.
Since also $\alpha\le\Leb(I_+)=\frac12$ and $\alpha\le1-t$, we obtain in both cases $\alpha\le\alpha_t$, where
\[
  \alpha_t\coloneqq
  \begin{cases}
    \min\{\frac12,\frac12-t+q\}, & t<\frac12,\\
    \min\{q,1-t\}, & t\ge\frac12.
  \end{cases}
\]

\emph{Step 3: the bathtub principle.}
On each half, the function $m_t$ has values in $[0,1]$ and a prescribed integral.
By the bathtub principle, see \cite[Thm.~1.14]{lieb2001analysis}, such a function is paired best with the nondecreasing function $f$ if it is the indicator of a right end segment of the interval.
Applying this on $I_+$ and on $I_-$ gives
\[
  \int_0^1f(x)m_t(x)\de x
  \le G(\alpha)\coloneqq\int_{1-\alpha}^1f(x)\de x+\int_{\alpha+t-\frac12}^{1/2}f(x)\de x.
\]
The function $G$ is Lipschitz, and its derivative $G'(\alpha)=f(1-\alpha)-f(\alpha+t-\frac12)$ is nonnegative almost everywhere because $1-\alpha\ge\frac12\ge\alpha+t-\frac12$.
Shifting mass from the lower half to the upper half is thus always profitable, and since $\alpha\le\alpha_t$, we obtain
\[
  \int_0^1f(x)m_t(x)\de x\le G(\alpha_t)=\int_{E_t}f(x)\de x,
  \qquad
  E_t\coloneqq[1-\alpha_t,1]\cup[\alpha_t+t-\tfrac12,\tfrac12).
\]

\emph{Step 4: identification of the optimal layers.}
Evaluating $\alpha_t$ in the four ranges of $t$ gives
\[
  E_t=
  \begin{cases}
    [t,1], & 0\le t<q,\\
    [q,\frac12)\cup[\frac12+t-q,1], & q\le t<\frac12,\\
    [t-\frac12+q,\frac12)\cup[1-q,1], & \frac12\le t<1-q,\\
    [t,1], & 1-q\le t<1.
  \end{cases}
\]
In each case, $E_t$ coincides up to a null set with $\{x:T^0_q(x)>t\}$, as one reads off from~\eqref{eq:T-maps}.
In other words, the coupling $\gamma^0_q$ attains the bound of Step~3 simultaneously in every layer.
Integrating over $t$ and using~\eqref{eq:layer-cake} for $\gamma^0_q$, we conclude that
\[
  \E_\gamma[f(X)Z]\le\int_0^1\int_0^1f(x)\1\{T^0_q(x)>t\}\de x\de t=\E_{\gamma^0_q}[f(X)Z].
\]

\emph{Lower bound.}
The random variable $Z'\coloneqq1-Z$ is uniform, and the law $\gamma'$ of $(X,Z')$ satisfies $\gamma'(Q)=\gamma(I_+\times I_-)=\frac12-q$.
The upper bound, applied to $\gamma'$ with $\frac12-q$ in place of $q$, gives $\E_{\gamma'}[f(X)Z']\le\E[f(X)T^0_{1/2-q}(X)]$.
Since $\E_\gamma[f(X)Z]=\E[f(X)]-\E_{\gamma'}[f(X)Z']$ and $T^1_q=1-T^0_{1/2-q}$, this is equivalent to the claimed lower bound.
\end{proof}

\begin{proof}[Proof of Theorem~\ref{thm:beta-eta}]
Let $C\in\CC$ with $\beta(C)=b$, put $q=\frac{1+b}{4}$, and let $\gamma$ be the associated coupling, so that $\gamma(Q)=q$.

\emph{Upper bound.}
By~\eqref{eq:coupling-moments}, $\nu(C)=12\,\E_\gamma[X^2Z]-2$, and Proposition~\ref{prop:constrained} with $f(x)=x^2$ bounds the expectation by its value under $\gamma^0_q$.
Integrating over the four linear pieces of $T^0_q$ gives
\begin{align*}
  \E_{\gamma^0_q}[X^2Z]
  &=\int_0^qx^3\de x+\int_q^{1/2}x^2\bigl(x+\tfrac12-q\bigr)\de x\\
  &\quad+\int_{1/2}^{1-q}x^2\bigl(x-\tfrac12+q\bigr)\de x+\int_{1-q}^1x^3\de x\\
  &=q^3-\tfrac32q^2+\tfrac34q+\tfrac18,
\end{align*}
and therefore
\[
  \nu(C)\le12\,\E_{\gamma^0_q}[X^2Z]-2=12q^3-18q^2+9q-\tfrac12=r_+(b),
\]
where the last identity follows from $q=\frac{1+b}{4}$ and~\eqref{eq:rpm-def}.
Since $\beta(C^\top)=\beta(C)$, the same bound applies to $\nu(C^\top)$, and averaging the two bounds gives $\eta(C)\le r_+(b)$ by~\eqref{eq:eta-def}.

\emph{Lower bound.}
Analogously, the lower bound of Proposition~\ref{prop:constrained} gives
\[
  \nu(C)\ge12\,\E_{\gamma^1_q}[X^2Z]-2=12q^3-1=r_-(b).
\]
The same bound holds for $\nu(C^\top)$, and hence $\eta(C)\ge r_-(b)$.

\emph{Attainment and fibers.}
The copulas $C^0_q$ and $C^1_q$ satisfy $\beta=b$ and are symmetric, so that $\eta(C^0_q)=\nu(C^0_q)=r_+(b)$ and $\eta(C^1_q)=\nu(C^1_q)=r_-(b)$.
The mixtures $(1-t)C^1_q+tC^0_q$, $t\in[0,1]$, have Blomqvist's beta $b$, and since $\eta$ is affine, their values of $\eta$ fill the fiber $[r_-(b),r_+(b)]$.
\end{proof}

The same argument with $f(x)=x$ recovers the bounds $r_-(b)\le\rho(C)\le r_+(b)$ of \cite{nelsen2001bounds}, attained by the same copulas $C^0_q$ and $C^1_q$.
Since these copulas are symmetric, $\nu$ and $\eta$ coincide on them, and by the computation above, $\rho$ takes the same values.
No pair of the three coefficients $\rho$, $\nu$ and $\eta$ can therefore separate the extremal copulas of the $(\beta,\cdot)$-regions.

\begin{corollary}[Deviation from Blomqvist's beta]
\label{cor:eta-beta-deviation}
Every copula $C$ satisfies $|\eta(C)-\beta(C)|\le\frac89$.
The bound is attained by $C^0_{1/6}$, which satisfies $\beta=-\frac13$ and $\eta=\frac59$, and by $C^1_{1/3}$, which satisfies $\beta=\frac13$ and $\eta=-\frac59$.
\end{corollary}

\begin{proof}
With $y\coloneqq1-b\in[0,2]$, the difference $r_+(b)-b=y-\frac{3}{16}y^3$ is maximal at $y=\frac43$, that is, at $b=-\frac13$ and $q=\frac16$, where it equals $\frac43-\frac49=\frac89$.
The difference $b-r_-(b)$ is treated in the same way, with maximum $\frac89$ at $b=\frac13$ and $q=\frac13$.
\end{proof}

\section{\texorpdfstring{The exact $(\xi,\eta)$-region}{The exact (xi,eta)-region}}
\label{sec:xi}

In contrast to Spearman's rho and Blomqvist's beta, Chatterjee's rank correlation is not an affine functional of the copula but a quadratic one, and it depends on the copula through its conditional distribution functions.
Our approach consists of three steps.
\begin{enumerate}[label=(\roman*),leftmargin=2em]
\item We write $\xi$ and $\eta$ in terms of the conditional distribution functions of the reflected coupling, in which $\eta$ is linear and $\xi$ is quadratic, see Lemma~\ref{lem:kappa}.
\item For a Lagrange multiplier $c>0$, we maximize $\eta-\frac c2\xi$.
A kernel of conditional distribution functions must be nondecreasing in the threshold $z$, and this is the only constraint that couples different thresholds.
Dropping it, the problem splits into independent convex problems, one for each threshold, whose solutions are clamped quadratic functions, see Lemma~\ref{lem:pointwise}.
\item We show that these solutions are nondecreasing in $z$ nevertheless, see Lemma~\ref{lem:families}.
They therefore define copulas, namely $C_c^\pm$, and the relaxation was exact.
\end{enumerate}
The remaining work consists in controlling the parametrization by $c$ and in computing the boundary curves in closed form.

\subsection{Conditional distribution functions}

We describe a coupling $\gamma\in\Gam$ by a jointly measurable kernel $\kappa:[0,1]^2\to[0,1]$ such that $\kappa(x,\cdot)$ is a version of the conditional distribution function $z\mapsto\gamma(Z\le z\mid X=x)$ for almost every $x$.
Since the first marginal of $\gamma$ is uniform, every such kernel satisfies
\begin{equation}
  \label{eq:kappa-normalization}
  \int_0^1\kappa(x,z)\de x=\gamma(Z\le z)=z,\qquad z\in[0,1].
\end{equation}
Conversely, let $\kappa(x,\cdot)$ be a distribution function on $[0,1]$ for almost every $x$ such that~\eqref{eq:kappa-normalization} holds.
Then $\gamma(\de x,\de z)\coloneqq\kappa(x,\de z)\de x$ defines a coupling in $\Gam$ with conditional distribution functions $\kappa$.
Throughout this section, we use the weight
\[
  a(x,z)\coloneqq x(x+2z)=(x+z)^2-z^2,\qquad(x,z)\in[0,1]^2.
\]

\begin{lemma}[Representation by conditional distribution functions]
\label{lem:kappa}
Let $C\in\CC$, and let $\kappa$ be a kernel of conditional distribution functions of the associated coupling.
Then
\begin{equation}
  \label{eq:kappa-rep}
  \begin{gathered}
  \xi(C)=6\int_0^1\!\int_0^1\kappa(x,z)^2\de x\de z-2,\\
  \eta(C)=3-6\int_0^1\!\int_0^1a(x,z)\,\kappa(x,z)\de x\de z.
  \end{gathered}
\end{equation}
Moreover, $\xi$ is convex on $\CC$ and lower semicontinuous with respect to uniform convergence.
\end{lemma}

\begin{proof}
We begin with $\xi$.
For $(U,V)\sim C$, the partial derivative $h\coloneqq\partial_1C$ exists almost everywhere with values in $[0,1]$, see \cite[Thm.~2.2.7]{nelsen2006introduction}, and $h(u,\cdot)$ is a version of the conditional distribution function of $V$ given $U=u$ for almost every $u$.
The conditional distribution of $V$ given $U=u$ has at most countably many atoms, so that
\[
  \kappa(x,z)=\Prob(1-V\le z\mid U=1-x)=1-h(1-x,1-z)
\]
for almost every $(x,z)$.
Since $\int_0^1h(u,v)\de u=v$, we have $\int_0^1\!\int_0^1h(u,v)\de u\de v=\frac12$, and expanding the square gives
\begin{align*}
  \int_0^1\!\int_0^1\kappa(x,z)^2\de x\de z
  &=1-2\int_0^1\!\int_0^1h(u,v)\de u\de v+\int_0^1\!\int_0^1h(u,v)^2\de u\de v\\
  &=\int_0^1\!\int_0^1h(u,v)^2\de u\de v,
\end{align*}
and~\eqref{eq:xi-def} yields the first identity in~\eqref{eq:kappa-rep}.

For $\eta$, fix $x\in[0,1]$ and put $g(z)\coloneqq x^2z+xz^2$.
Since $g(0)=0$ and $g'(t)=a(x,t)$, we can write $g(Z)=\int_0^1a(x,t)\1\{Z>t\}\de t$, and taking conditional expectations gives
\[
  \E\bigl[XZ(X+Z)\bigm|X=x\bigr]=\int_0^1a(x,t)\bigl(1-\kappa(x,t)\bigr)\de t.
\]
Integrating over $x$ and using $\int_0^1\!\int_0^1a(x,t)\de x\de t=\frac13+\frac12=\frac56$, we obtain
\[
  \E\bigl[XZ(X+Z)\bigr]=\frac56-\int_0^1\!\int_0^1a(x,t)\,\kappa(x,t)\de x\de t,
\]
and~\eqref{eq:coupling-moments} yields the second identity.

The map $C\mapsto\partial_1C$ is affine, so $\xi$ is convex by~\eqref{eq:xi-def}.
For the lower semicontinuity, let $C_n\to C$ uniformly, and let $\phi$ be smooth with compact support in $(0,1)^2$.
Integration by parts in the first variable gives
\[
  \int_0^1\!\int_0^1\partial_1C_n(u,v)\,\phi(u,v)\de u\de v
  =-\int_0^1\!\int_0^1C_n(u,v)\,\partial_1\phi(u,v)\de u\de v.
\]
As $n\to\infty$, the right-hand side converges to the corresponding expression for $C$, which in turn equals the left-hand side for $C$.
Since $|\partial_1C_n|\le1$ and such functions $\phi$ are dense in $L^2([0,1]^2)$, it follows that $\partial_1C_n\to\partial_1C$ weakly in $L^2([0,1]^2)$.
The weak lower semicontinuity of the norm then yields $\xi(C)\le\liminf_n\xi(C_n)$.
\end{proof}

\subsection{The relaxed problem}

For a Lagrange multiplier $c>0$, Lemma~\ref{lem:kappa} expresses the Lagrangian of the upper boundary problem as
\begin{equation}
  \label{eq:lagrange-identity}
  \eta(C)-\frac c2\,\xi(C)=3+c-6\int_0^1\Phi_z\bigl(\kappa(\cdot,z)\bigr)\de z,
\end{equation}
where, for each threshold $z\in[0,1]$ and each square-integrable $f:[0,1]\to\R$,
\begin{equation}
  \label{eq:Phi-z}
  \Phi_z(f)\coloneqq\int_0^1\Bigl(a(x,z)f(x)+\frac c2f(x)^2\Bigr)\de x.
\end{equation}
Maximizing the Lagrangian thus means minimizing $\int_0^1\Phi_z(\kappa(\cdot,z))\de z$ over all kernels of copulas.
Such a kernel is subject to three constraints: it takes values in $[0,1]$, it satisfies the normalization~\eqref{eq:kappa-normalization} at every threshold, and it is nondecreasing in $z$.
Only the last constraint links different thresholds.
If we drop it, the minimization decouples into one problem for each $z$, in which $\Phi_z$ is minimized over functions with values in $[0,1]$ and integral $z$.
Each of these problems is strictly convex, and its solution is obtained by clamping the unconstrained minimizer, a quadratic function of the rank sum $x+z$, to the interval $[0,1]$.

\begin{lemma}[Level-wise minimizers]
\label{lem:pointwise}
Let $c>0$ and $z\in(0,1)$.
There are unique numbers $B=B_c(z)>z$ and $L=L_c(z)\in\R$ such that the functions
\[
  \kappa^+(x)\coloneqq\clip{\frac{B^2-(x+z)^2}{c}},
  \qquad
  \kappa^-(x)\coloneqq\clip{\frac{(x+z)^2-L}{c}},
  \qquad x\in[0,1],
\]
satisfy $\int_0^1\kappa^\pm(x)\de x=z$.
For every measurable $f:[0,1]\to[0,1]$ with $\int_0^1f(x)\de x=z$,
\begin{equation}
  \label{eq:pointwise}
  \begin{gathered}
  \Phi_z(f)\ge\Phi_z(\kappa^+),\\
  \int_0^1\Bigl(\frac c2f(x)^2-a(x,z)f(x)\Bigr)\de x\ge\int_0^1\Bigl(\frac c2\kappa^-(x)^2-a(x,z)\kappa^-(x)\Bigr)\de x,
  \end{gathered}
\end{equation}
with equality if and only if $f=\kappa^+$ or $f=\kappa^-$ almost everywhere, respectively.
\end{lemma}

\begin{proof}
We first show that $B$ exists and is unique.
The map $B'\mapsto\int_0^1\clip{(B'^2-(x+z)^2)/c}\de x$ is continuous and nondecreasing on $[z,\infty)$, it equals $0$ at $B'=z$, and it equals $1$ for large $B'$.
Suppose that its value at some $B'$ lies in $(0,1)$.
Then the integrand, which is continuous and nonincreasing in $x$, is neither identically $0$ nor identically $1$, so it takes values in $(0,1)$ on an interval of positive length.
On this interval, the integrand strictly increases with $B'$, so the map is strictly increasing wherever its value lies in $(0,1)$.
Hence there is exactly one $B>z$ with $\int_0^1\kappa^+(x)\de x=z$.
The same argument applies to $L$.

Now let $f$ be as in the statement and put $\lambda\coloneqq B^2-z^2$.
Expanding the square and using $\int_0^1f(x)\de x=z=\int_0^1\kappa^+(x)\de x$, we obtain
\begin{align*}
  \Phi_z(f)-\Phi_z(\kappa^+)
  &=\int_0^1\bigl(a(x,z)+c\kappa^+(x)\bigr)\bigl(f(x)-\kappa^+(x)\bigr)\de x\\
  &\quad+\frac c2\int_0^1\bigl(f(x)-\kappa^+(x)\bigr)^2\de x\\
  &=\int_0^1\bigl(a(x,z)+c\kappa^+(x)-\lambda\bigr)\bigl(f(x)-\kappa^+(x)\bigr)\de x\\
  &\quad+\frac c2\int_0^1\bigl(f(x)-\kappa^+(x)\bigr)^2\de x.
\end{align*}
The first factor of the first integrand equals $(x+z)^2-B^2+c\kappa^+(x)$, and we check its sign against that of the second factor.
Where $0<\kappa^+(x)<1$, the first factor vanishes.
Where $\kappa^+(x)=0$, it is nonnegative, and so is $f(x)-\kappa^+(x)=f(x)$.
Where $\kappa^+(x)=1$, it is nonpositive, and so is $f(x)-1$.
The first integrand is therefore nonnegative, which proves the first inequality in~\eqref{eq:pointwise}, with equality if and only if $f=\kappa^+$ almost everywhere.
The second inequality follows in the same way with $-a$ in place of $a$ and with $\lambda\coloneqq z^2-L$.
\end{proof}

\subsection{The extremal families}

For $c>0$, we collect the level-wise minimizers in the kernels
\begin{equation}
  \label{eq:kappa-pm}
  \kappa_c^+(x,z)\coloneqq\clip{\frac{B_c(z)^2-(x+z)^2}{c}},
  \qquad
  \kappa_c^-(x,z)\coloneqq\clip{\frac{(x+z)^2-L_c(z)}{c}}
\end{equation}
for $(x,z)\in[0,1]\times(0,1)$, and we put $\kappa_c^\pm(x,0)\coloneqq0$ and $\kappa_c^\pm(x,1)\coloneqq1$.
The key step of the proof is that these kernels are nondecreasing in $z$, so that dropping the monotonicity constraint was harmless.

\begin{lemma}[The extremal families]
\label{lem:families}
Let $c>0$.
For every $x\in(0,1)$, the functions $z\mapsto\kappa_c^\pm(x,z)$ are continuous and nondecreasing distribution functions on $[0,1]$.
Consequently, $\kappa_c^+$ and $\kappa_c^-$ are the conditional distribution functions of couplings $\gamma_c^+,\gamma_c^-\in\Gam$.
The copulas $C_c^+$ and $C_c^-$ associated with the kernels~\eqref{eq:kappa-pm} are absolutely continuous and satisfy~\eqref{eq:families} with $a_c(v)=B_c(1-v)^2-c$ and $b_c(v)=L_c(1-v)+c$, which are the unique solutions of~\eqref{eq:ab-def}.
\end{lemma}

\begin{proof}
We treat $\kappa_c^+$ in detail and write $B=B_c$.

\emph{Differentiability of $B$.}
For $z\in(0,1)$ and $\tilde B\ge z$, the substitution $y=x+z$ gives
\[
  \begin{gathered}
  \int_0^1\clip{\frac{\tilde B^2-(x+z)^2}{c}}\de x=I(\tilde B,z)\coloneqq\int_z^{1+z}k_{\tilde B}(y)\de y,\\
  k_{\tilde B}(y)\coloneqq\clip{\frac{\tilde B^2-y^2}{c}}.
  \end{gathered}
\]
The variable $y$ is the rank sum, and $k_{\tilde B}(y)\in(0,1)$ precisely for $y$ in an interval, which we call the active interval.
The function $I$ is continuously differentiable, with
\[
  \partial_zI(\tilde B,z)=k_{\tilde B}(1+z)-k_{\tilde B}(z)
  \qquad\text{and}\qquad
  \partial_{\tilde B}I(\tilde B,z)=\frac{2\tilde B}{c}\,\ell(\tilde B,z),
\]
where $\ell(\tilde B,z)$ denotes the length of the active interval within $[z,1+z]$.
By the proof of Lemma~\ref{lem:pointwise}, $\ell(B(z),z)>0$.
Since $B(z)$ solves $I(B(z),z)=z$, the implicit function theorem shows that $B$ is continuously differentiable on $(0,1)$ with
\begin{equation}
  \label{eq:B-prime}
  B(z)B'(z)=\frac{c\,\bigl(1-k(1+z)+k(z)\bigr)}{2\ell},
\end{equation}
where $k\coloneqq k_{B(z)}$ and $\ell\coloneqq\ell(B(z),z)$.

\emph{The key inequality.}
Put $A\coloneqq\sqrt{\max\{B^2-c,0\}}$.
Then $k(y)\in(0,1)$ precisely for $y\in(A,B)$, and the active interval within $[z,1+z]$ is $(\max\{A,z\},\min\{B,1+z\})$.
We claim that
\begin{equation}
  \label{eq:key-inequality}
  B(z)B'(z)\ge\min\{B(z),1+z\},\qquad z\in(0,1).
\end{equation}
We distinguish four cases according to the position of the active interval, using that $B^2-A^2\le c$, with equality if $A>0$.
\begin{enumerate}[label=(\roman*),leftmargin=2em]
\item \emph{$A\ge z$ and $B\le1+z$.}
Then $k(z)=1$, $k(1+z)=0$ and $\ell=B-A$.
Since $c=B^2-A^2$, formula~\eqref{eq:B-prime} gives $BB'=A+B\ge B$.
\item \emph{$A<z$ and $B\le1+z$.}
Then $k(z)=(B^2-z^2)/c$, $k(1+z)=0$ and $\ell=B-z$.
Hence $BB'\ge B$ is equivalent to $c\ge(B-z)^2$.
This holds because $c>B^2-z^2\ge(B-z)^2$.
\item \emph{$A\ge z$ and $B>1+z$.}
Then $k(z)=1$, $k(1+z)=(B^2-(1+z)^2)/c$ and $\ell=1+z-A$.
Hence $BB'\ge1+z$ is equivalent to $c\ge(1+z-A)^2$.
This holds because $c=B^2-A^2$ exceeds $(1+z)^2-A^2\ge(1+z-A)^2$.
\item \emph{$A<z$ and $B>1+z$.}
Then $\ell=1$ and $BB'=\frac12\bigl(c+(1+z)^2-z^2\bigr)$.
Hence $BB'\ge1+z$ is equivalent to $c\ge1$.
This holds because $c>B^2-z^2>(1+z)^2-z^2\ge1$.
\end{enumerate}

\emph{Monotonicity in $z$.}
Fix $x\in(0,1)$ and put $g(z)\coloneqq B(z)^2-(x+z)^2$.
Then $\kappa_c^+(x,z)=\clip{g(z)/c}$.
Suppose that $0<g(z)<c$.
Then the rank sum $x+z$ lies in the active interval, and hence $x+z<\min\{B(z),1+z\}$.
By~\eqref{eq:key-inequality}, we then have
\[
  g'(z)=2B(z)B'(z)-2(x+z)\ge2\min\{B(z),1+z\}-2(x+z)>0.
\]
The function $z\mapsto\clip{g(z)/c}$ is locally Lipschitz on $(0,1)$.
Its derivative equals $g'/c\ge0$ on the open set $\{0<g<c\}$ and vanishes almost everywhere on the complement of this set.
Hence $\kappa_c^+(x,\cdot)$ is nondecreasing on $(0,1)$.

\emph{Boundary values.}
The function $\kappa_c^+(\cdot,z)$ is nonincreasing in $x$ and has integral $z$.
As $z\downarrow0$, the monotone limit of $\kappa_c^+(x,z)$ is therefore a nonincreasing function of $x$ with integral $0$, which vanishes for every $x\in(0,1)$.
In the same way, the limit as $z\uparrow1$ equals $1$ for every $x\in(0,1)$.
Hence $\kappa_c^+(x,\cdot)$ is a continuous distribution function on $[0,1]$ for every $x\in(0,1)$.

\emph{The lower family.}
For $\kappa_c^-$, we use $I(\tilde L,z)\coloneqq\int_z^{1+z}\clip{(y^2-\tilde L)/c}\de y$ instead, for which $\partial_{\tilde L}I=-\ell/c$.
With $k(y)\coloneqq\clip{(y^2-L(z))/c}$ and the length $\ell$ of its active interval, the same computation gives
\[
  L'(z)=-\frac{c\,\bigl(1-k(1+z)+k(z)\bigr)}{\ell}\le0.
\]
Hence $\partial_z\bigl((x+z)^2-L(z)\bigr)\ge2(x+z)>0$, and $\kappa_c^-(x,\cdot)$ is nondecreasing by the same argument as before.
The boundary values follow as for $\kappa_c^+$.

\emph{Conclusion.}
The functions $\kappa_c^\pm(x,\cdot)$ are monotone, continuous on $[0,1]$ and locally Lipschitz on $(0,1)$, hence absolutely continuous.
By the remarks at the beginning of this section, $\kappa_c^+$ and $\kappa_c^-$ therefore define couplings in $\Gam$, and the associated copulas $C_c^\pm$ are absolutely continuous.
The proof of Lemma~\ref{lem:kappa} shows that $\partial_1C(u,v)=1-\kappa(1-u,1-v)$ for almost all $(u,v)$.
Since $1-\clip{t}=\clip{1-t}$, the copulas satisfy~\eqref{eq:families} with the stated $a_c$ and $b_c$.
Substituting $x=1-u$ and $z=1-v$ in the normalization of Lemma~\ref{lem:pointwise} shows that these numbers solve~\eqref{eq:ab-def}, and they are unique there by the uniqueness of $B_c(z)$ and $L_c(z)$.
\end{proof}

\begin{figure}[t!]
  \centering
  \includegraphics[width=0.82\textwidth]{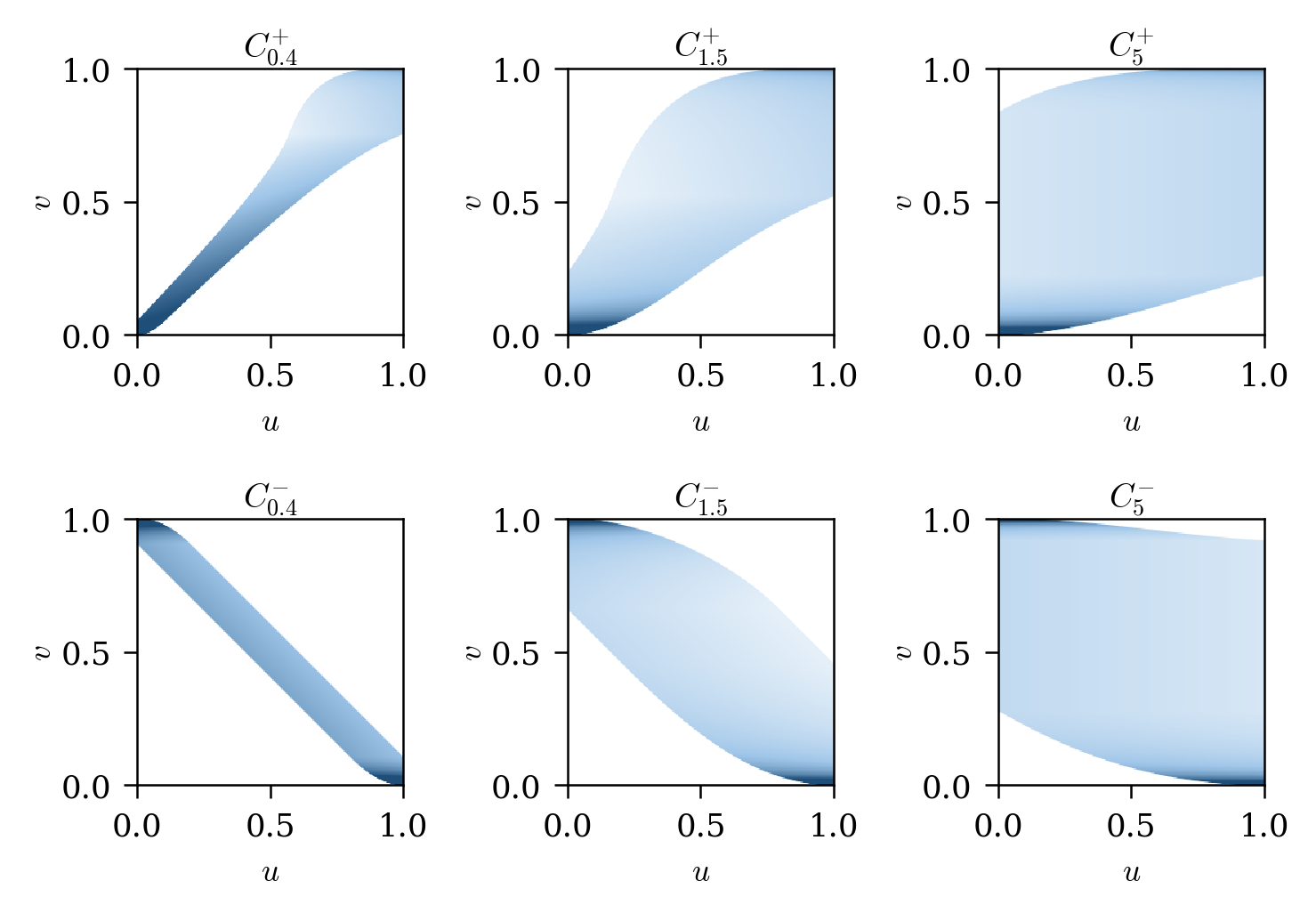}
  \caption{Densities of the extremal copulas $C_c^+$ (top) and $C_c^-$ (bottom) of Theorem~\ref{thm:xi-eta} for $c=0.4$, $1.5$ and $5$, computed from~\eqref{eq:families} by solving~\eqref{eq:ab-def} numerically, darker shades indicating larger values and white indicating zero density.
  The corresponding boundary points $(\xi,\eta)$ are approximately $(0.77,0.98)$, $(0.41,0.82)$ and $(0.07,0.37)$ for $C_c^+$, and $(0.81,-0.98)$, $(0.41,-0.80)$ and $(0.07,-0.36)$ for $C_c^-$.
  As $c$ decreases, the upper copulas concentrate near the main diagonal and the lower ones near the secondary diagonal.}
  \label{fig:xi-extremals}
\end{figure}

The extremal copulas have a simple dependence structure.
Since $(2-u-v)^2$ decreases in $u$, the conditional distribution function $\partial_1C_c^+(u,v)$ is nonincreasing in $u$, so that $V$ is stochastically increasing in $U$ under $C_c^+$.
In the same way, $V$ is stochastically decreasing in $U$ under $C_c^-$.
Figure~\ref{fig:xi-extremals} shows the densities of both families for three values of $c$.

\subsection{A square-root bound}

If we also drop the constraint $0\le\kappa\le1$ in the relaxed problem, each level-wise problem becomes an unconstrained quadratic problem with an explicit value.
This yields a universal inequality between $\eta$ and $\xi$, which is sharp near independence and which we use in the proof of Theorem~\ref{thm:xi-eta}.

\begin{proposition}[A sharp square-root bound]
\label{prop:sqrt}
For every copula $C$,
\[
  |\eta(C)|\le\sqrt{\tfrac{11}{5}\,\xi(C)},
\]
with equality if and only if $C=\Pi$.
Moreover, $\sup\{\pm\eta(C):C\in\CC,\ \xi(C)=x\}=\sqrt{11x/5}\,(1+o(1))$ as $x\downarrow0$, so the constant $\frac{11}{5}$ is best possible.
\end{proposition}

\begin{proof}
\emph{The inequality.}
Fix $c>0$ and $z\in[0,1]$, and let $f\in L^2([0,1])$ satisfy $\int_0^1f(x)\de x=z$.
For a uniform random variable $X$, the mean and the variance of $a(X,z)$ are
\begin{gather*}
  \bar a(z)\coloneqq\int_0^1a(x,z)\de x=\frac13+z,\\
  V(z)\coloneqq\int_0^1a(x,z)^2\de x-\bar a(z)^2=\frac{4}{45}+\frac z3+\frac{z^2}{3}.
\end{gather*}
Completing the square in~\eqref{eq:Phi-z} gives
\[
  \Phi_z(f)=\frac c2\int_0^1\Bigl(f(x)+\frac{a(x,z)}{c}\Bigr)^2\de x-\frac{1}{2c}\int_0^1a(x,z)^2\de x.
\]
The function $f+a(\cdot,z)/c$ has integral $z+\bar a(z)/c$, so Jensen's inequality yields
\[
  \Phi_z(f)\ge\frac c2\Bigl(z+\frac{\bar a(z)}{c}\Bigr)^2-\frac{1}{2c}\int_0^1a(x,z)^2\de x
  =\frac{c}{2}z^2+z\,\bar a(z)-\frac{V(z)}{2c}.
\]
We apply this bound to $f=\kappa(\cdot,z)$, where $\kappa$ is the kernel of a copula $C$, and integrate over $z$.
Since $\int_0^1z\,\bar a(z)\de z=\frac12$ and $\int_0^1V(z)\de z=\frac{11}{30}$, we obtain
\[
  \int_0^1\Phi_z\bigl(\kappa(\cdot,z)\bigr)\de z\ge\frac c6+\frac12-\frac{11}{60c}.
\]
By~\eqref{eq:lagrange-identity}, this means
\begin{equation}
  \label{eq:relaxed}
  \eta(C)-\frac c2\,\xi(C)\le\frac{11}{10\,c}.
\end{equation}
The same computation with $-a$ in place of $a$ gives $-\eta(C)-\frac c2\,\xi(C)\le\frac{11}{10c}$.

Both inequalities are strict.
Indeed, equality in~\eqref{eq:relaxed} would require equality in Jensen's inequality for almost every $z$, that is, $\kappa(x,z)=z+(\bar a(z)-a(x,z))/c$ for almost all $(x,z)$.
Since $a(1,z)-\bar a(z)=\frac23+z$, this function is negative for $x$ close to $1$ and $z$ close to $0$, which is impossible for a kernel of conditional distribution functions.
The same argument applies to the second inequality.

If $\xi(C)>0$, we choose $c=\sqrt{11/(5\xi(C))}$, for which $\frac{11}{10c}+\frac c2\xi(C)=\sqrt{11\xi(C)/5}$, and the two strict inequalities become $|\eta(C)|<\sqrt{11\xi(C)/5}$.
If $\xi(C)=0$, then $C=\Pi$ by \cite[Thm.~1.1]{chatterjee2020}, and both sides vanish.

\emph{Sharpness.}
We perturb the independence kernel $\kappa(x,z)=z$ in the direction suggested by the proof.
Let $\delta\in(0,\frac12)$, put $\psi_\delta(z)\coloneqq\min\{1,z/\delta,(1-z)/\delta\}$, and consider
\[
  \kappa_\varepsilon(x,z)\coloneqq z+\varepsilon\,\psi_\delta(z)\bigl(\bar a(z)-a(x,z)\bigr),
  \qquad|\varepsilon|\le\frac\delta4.
\]
We check that $\kappa_\varepsilon$ is the kernel of a coupling.
First, $\int_0^1\kappa_\varepsilon(x,z)\de x=z$ because $a(\cdot,z)-\bar a(z)$ has mean zero.
Second, the bounds $|\bar a(z)-a(x,z)|\le\frac53$ and $\psi_\delta(z)\le\min\{z,1-z\}/\delta$ imply
\[
  \bigl(1-\tfrac{5}{12}\bigr)z\le\kappa_\varepsilon(x,z)\le z+\tfrac{5}{12}(1-z),
\]
so $\kappa_\varepsilon$ takes values in $[0,1]$.
Third, $|\psi_\delta'|\le1/\delta$ and $|\partial_z(\bar a(z)-a(x,z))|=|1-2x|\le1$ give
\[
  \partial_z\kappa_\varepsilon(x,z)\ge1-\frac{5}{12}-\frac18>0,
\]
so $\kappa_\varepsilon(x,\cdot)$ is nondecreasing.
Hence $\kappa_\varepsilon$ is the kernel of a copula $C_\varepsilon$.
Finally, we have
\[
  \int_0^1\bigl(\bar a(z)-a(x,z)\bigr)\,a(x,z)\de x=-V(z),
\]
and Lemma~\ref{lem:kappa} gives
\[
  \xi(C_\varepsilon)=6\varepsilon^2\int_0^1\psi_\delta(z)^2V(z)\de z,
  \qquad
  \eta(C_\varepsilon)=6\varepsilon\int_0^1\psi_\delta(z)V(z)\de z.
\]
For every sufficiently small $x>0$, the choice $\varepsilon=\pm\bigl(x/(6\int_0^1\psi_\delta^2V\de z)\bigr)^{1/2}$ is admissible and yields copulas with $\xi=x$ and
\[
  \pm\eta=\sqrt{x}\cdot\frac{\sqrt6\int_0^1\psi_\delta(z)V(z)\de z}{\bigl(\int_0^1\psi_\delta(z)^2V(z)\de z\bigr)^{1/2}}.
\]
As $\delta\downarrow0$, the fraction converges to $\bigl(6\int_0^1V(z)\de z\bigr)^{1/2}=\sqrt{11/5}$.
Together with the inequality, this proves the asymptotic statement.
\end{proof}

\subsection{The boundary curves in closed form}

The level-wise normalization $\int_0^1\kappa_c^\pm(x,z)\de x=z$ can be solved explicitly in each of the four configurations of the active interval that occur in the proof of Lemma~\ref{lem:families}.
Changing the variable of integration from $z$ to a suitable parameter of each configuration, the integrals in~\eqref{eq:kappa-rep} become integrals of rational functions and can be evaluated in closed form.
We present the resulting formulas here and give the configurations and the derivation in Appendix~\ref{app:closed}.
The formulas change their form at the branch points
\[
  c_+\coloneqq1+\frac{2}{\sqrt3}\approx2.155,
  \qquad
  c_-\coloneqq1+\sqrt{\frac23}\approx1.816.
\]

For $c\ge c_+$, let
\begin{equation}
  \label{eq:up-large}
  \begin{aligned}
  x_+(c)&\coloneqq\frac{2760c^4-1080c^3-1921c^2+1}{405c^2(c^2-1)}-\frac{8(c-1)}{\sqrt c}\arctan\frac{1}{\sqrt c}\\
  &\quad+\frac{16c}{27}\log\frac{c+1}{c-1}-8\log\frac{c+1}{c},\\
  e_+(c)&\coloneqq\frac{420c^4+1440c^3-c^2-1980c+1}{405c(c^2-1)}-\frac{4\sqrt c\,(c+3)}{3}\arctan\frac{1}{\sqrt c}\\
  &\quad+\frac{4c^2}{27}\log\frac{c+1}{c-1}.
  \end{aligned}
\end{equation}
For $0<c\le c_+$, let $P\ge\sqrt{c_+}$ be the unique solution of
\begin{equation}
  \label{eq:P-def}
  \sqrt c\,(3P^4+6P^2-1)=12P^3,
\end{equation}
put $K\coloneqq3P^4+6P^2-1$, and let
\begin{equation}
  \label{eq:up-small}
  \begin{aligned}
  x_+(c)&\coloneqq\frac{4(108P^8+277P^6+360P^4-144P^2+12)}{3K^2}\log P\\
  &\quad-\frac{8(27P^8+76P^6+90P^4-36P^2+3)}{3K^2}\log(1+P^2)\\
  &\quad+\frac{2(3P^4-12P^3+6P^2-1)(3P^4+12P^3+6P^2-1)}{3P^3K}\arctan\frac1P\\
  &\quad-\frac{810P^{10}-26460P^8-41460P^6+7068P^4+346P^2-21}{90P^2K^2}+k_\xi c,\\
  e_+(c)&\coloneqq\frac{3072P^{12}}{K^4}\log(1+P^2)-\frac{7440P^{12}}{K^4}\log P\\
  &\quad-\frac{48P^3(9P^8+84P^6+30P^4-12P^2+1)}{K^3}\arctan\frac1P-k_\eta c^2\\
  &\quad+\frac{405P^{16}+9720P^{14}+81540P^{12}+103800P^{10}}{5K^4}\\
  &\quad-\frac{8574P^8+17048P^6-5178P^4+600P^2-25}{5K^4},
  \end{aligned}
\end{equation}
with the constants
\begin{align*}
  k_\xi&\coloneqq\frac{559}{288}-\frac{157\sqrt3}{108}-\frac{32}{27}\log(\sqrt3-1)+\frac18\log\Bigl(1+\frac{2}{\sqrt3}\Bigr)\approx-0.1113,\\
  k_\eta&\coloneqq-\frac{565}{1152}+\frac{157\sqrt3}{432}+\frac{8}{27}\log(\sqrt3-1)-\frac{1}{32}\log\Bigl(1+\frac{2}{\sqrt3}\Bigr)\approx0.0226.
\end{align*}
For $c\ge c_-$, with $w\coloneqq\sqrt c$, let
\begin{equation}
  \label{eq:lo-large}
  \begin{aligned}
  x_-(c)&\coloneqq\frac{480w^8-1890w^6+1381w^4-1}{405w^4(w^4-1)}-\frac{4w^2-1}{w}\arctan\frac1w-8\log w\\
  &\quad+\frac{(2w+3)^2(8w+3)}{54w}\log(w+1)+\frac{(2w-3)^2(8w-3)}{54w}\log(w-1)\\
  &\quad-\frac{2(8w^2-27)}{27}\log(w^2+1),\\
  e_-(c)&\coloneqq\frac{w(4w^2+3)}{6}\arctan\frac1w-\frac{120w^8+225w^6+811w^4-1170w^2-1}{405w^2(w^4-1)}\\
  &\quad-\frac{w(2w+3)^2(4w-3)}{108}\log(w+1)-\frac{w(2w-3)^2(4w+3)}{108}\log(w-1)\\
  &\quad+\frac{4w^4}{27}\log(w^2+1).
  \end{aligned}
\end{equation}
For $0<c\le c_-$, let $p\ge\sqrt{c_-}$ be the unique solution of
\begin{equation}
  \label{eq:p-def}
  \sqrt c\,(3p^4+1)=6p^3,
\end{equation}
put $R\coloneqq3p^4+1$, and let
\begin{equation}
  \label{eq:lo-small}
  \begin{aligned}
  x_-(c)&\coloneqq\frac{(p+1)^4(3p^2-2p+1)^2(3p^4+16p^3+1)}{12p^3R^2}\log(p+1)\\
  &\quad-\frac{(p-1)^4(3p^2+2p+1)^2(3p^4-16p^3+1)}{12p^3R^2}\log(p-1)\\
  &\quad+\frac{2(27p^8-32p^6+18p^4+3)}{3R^2}\log(p^2+1)-8\log p\\
  &\quad+\frac{(3p^4-12p^3+1)(3p^4+12p^3+1)}{6p^3R}\arctan\frac1p\\
  &\quad-\frac{2(405p^8-420p^6+54p^4+220p^2-27)}{45p^2R^2},\\
  e_-(c)&\coloneqq\frac{3p^3(p+1)^4(3p^2-2p+1)^2(3p^4-8p^3+1)}{2R^4}\log(p+1)\\
  &\quad-\frac{3p^3(p-1)^4(3p^2+2p+1)^2(3p^4+8p^3+1)}{2R^4}\log(p-1)\\
  &\quad+\frac{192p^{12}}{R^4}\log(p^2+1)\\
  &\quad+\frac{3p^3(3p^2+1)(3p^6+15p^4-3p^2+1)}{R^3}\arctan\frac1p\\
  &\quad-\frac{405p^{16}+810p^{14}-270p^{12}+2100p^{10}}{5R^4}\\
  &\quad-\frac{108p^8-110p^6+126p^4+15}{5R^4}.
  \end{aligned}
\end{equation}

\begin{proposition}[Closed forms]
\label{prop:closed}
For every $c>0$,
\begin{gather*}
  \xi(C_c^+)=x_+(c),\qquad\eta(C_c^+)=e_+(c),\\
  \xi(C_c^-)=x_-(c),\qquad\eta(C_c^-)=e_-(c).
\end{gather*}
The two branches of each function agree at $c=c_\pm$, and as $c\to\infty$,
\begin{align*}
  x_+(c)&=\frac{11}{5c^2}-\frac{272}{105c^3}+O(c^{-4}),
  &e_+(c)&=\frac{11}{5c}-\frac{68}{35c^2}+O(c^{-3}),\\
  x_-(c)&=\frac{11}{5c^2}-\frac{328}{105c^3}+O(c^{-4}),
  &e_-(c)&=-\frac{11}{5c}+\frac{82}{35c^2}+O(c^{-3}).
\end{align*}
\end{proposition}

The proof is given in the appendix.
The expansions show that the fibers of the $(\xi,\eta)$-region are not symmetric about zero.
Indeed, eliminating $c$ gives
\begin{equation}
  \label{eq:eta-expansion}
  \eta_+(x)=\sqrt{\tfrac{11}{5}x}-\tfrac{68}{231}x+O(x^{3/2}),
  \qquad
  \eta_-(x)=-\sqrt{\tfrac{11}{5}x}+\tfrac{82}{231}x+O(x^{3/2})
\end{equation}
as $x\downarrow0$ for the boundary functions of Theorem~\ref{thm:xi-eta}.

\subsection{Proof of Theorem~\ref{thm:xi-eta}}

By Proposition~\ref{prop:closed}, the curves $c\mapsto(x_\pm(c),e_\pm(c))$ of the theorem are the curves $c\mapsto(\xi(C_c^\pm),\eta(C_c^\pm))$ traced by the extremal families.
We prove the theorem for the upper boundary and write $C_c\coloneqq C_c^+$.
The lower boundary is treated analogously, with the corresponding changes indicated in each step.

\emph{Step 1: a Lagrangian inequality.}
We show that for every $C\in\CC$ and $c>0$,
\begin{equation}
  \label{eq:lagrange}
  \eta(C)-\frac c2\,\xi(C)\le\eta(C_c)-\frac c2\,\xi(C_c),
\end{equation}
with equality if and only if $C=C_c$.
Let $\kappa$ be the kernel of $C$.
For almost every $z$, the function $\kappa(\cdot,z)$ takes values in $[0,1]$ and has integral $z$, so Lemma~\ref{lem:pointwise} gives $\Phi_z(\kappa(\cdot,z))\ge\Phi_z(\kappa_c^+(\cdot,z))$.
Integrating over $z$ and using~\eqref{eq:lagrange-identity} for both $C$ and $C_c$, which is a copula by Lemma~\ref{lem:families}, we obtain~\eqref{eq:lagrange}.
Equality forces $\kappa(\cdot,z)=\kappa_c^+(\cdot,z)$ almost everywhere for almost every $z$, that is, $C=C_c$.
For the lower boundary, the corresponding inequality reads $\eta(C)+\frac c2\xi(C)\ge\eta(C_c^-)+\frac c2\xi(C_c^-)$.

\emph{Step 2: monotonicity of the parametrization.}
Let $c_1<c_2$.
The copulas $C_{c_1}$ and $C_{c_2}$ differ, because on the nondegenerate active intervals their kernels have the different slopes $-2(x+z)/c_1$ and $-2(x+z)/c_2$ in $x$.
Applying~\eqref{eq:lagrange} with strict inequality to $C=C_{c_2}$ and $c=c_1$, and to $C=C_{c_1}$ and $c=c_2$, gives
\begin{gather*}
  \eta(C_{c_2})-\tfrac{c_1}{2}\xi(C_{c_2})<\eta(C_{c_1})-\tfrac{c_1}{2}\xi(C_{c_1}),\\
  \eta(C_{c_1})-\tfrac{c_2}{2}\xi(C_{c_1})<\eta(C_{c_2})-\tfrac{c_2}{2}\xi(C_{c_2}).
\end{gather*}
Adding the two inequalities yields $(c_2-c_1)\bigl(\xi(C_{c_1})-\xi(C_{c_2})\bigr)>0$, so $c\mapsto\xi(C_c)$ is strictly decreasing.
The first of the two inequalities then shows that $c\mapsto\eta(C_c)$ is strictly decreasing as well.

\emph{Step 3: continuity and range of the parametrization.}
For fixed $z\in(0,1)$, the implicit function theorem, applied to the equation $I(B,z)=z$ from the proof of Lemma~\ref{lem:families} with $c$ as an additional variable, shows that $B_c(z)$ depends continuously on $c$.
Hence $\kappa_c^+$ depends pointwise continuously on $c$, and dominated convergence in~\eqref{eq:kappa-rep} shows that $\xi(C_c)$ and $\eta(C_c)$ are continuous in $c$.
Inequality~\eqref{eq:lagrange} with $C=\Pi$ gives $\frac c2\xi(C_c)\le\eta(C_c)\le1$, so $\xi(C_c)\to0$ as $c\to\infty$.
With $C=M$, it gives $\eta(C_c)\ge1-\frac c2$, so $\eta(C_c)\to1$ as $c\downarrow0$.
Suppose that $\xi(C_{c_n})\le x_0<1$ along a sequence $c_n\downarrow0$.
Since $\CC$ is compact, a subsequence converges uniformly to some $\widetilde C\in\CC$, which satisfies $\xi(\widetilde C)\le x_0$ by Lemma~\ref{lem:kappa} and $\eta(\widetilde C)=1$ by continuity of the affine functional $\eta$.
But
\[
  \eta(M)-\eta(\widetilde C)=12\int_0^1\!\int_0^1(2-u-v)\bigl(M(u,v)-\widetilde C(u,v)\bigr)\de u\de v
\]
with $M\ge\widetilde C$ shows that $\eta(\widetilde C)=1$ forces $\widetilde C=M$, contradicting $\xi(M)=1$.
Therefore $\xi(C_c)\to1$ as $c\downarrow0$, and $c\mapsto\xi(C_c)=x_+(c)$ is a continuous and strictly decreasing bijection from $(0,\infty)$ onto $(0,1)$.
In particular, $\eta_+$ is well defined, and by Step~2, it is strictly increasing on $(0,1)$.
For the lower family, one uses $W$ in place of $M$ and the bound $\eta(C_c^-)\le-1+\frac c2$.

\emph{Step 4: the bounds.}
Let $x\in(0,1)$, let $c$ be such that $\xi(C_c)=x$, and let $C\in\CC$ with $\xi(C)\le x$.
Then~\eqref{eq:lagrange} gives
\[
  \eta(C)\le\eta(C_c)-\frac c2\bigl(x-\xi(C)\bigr)\le\eta_+(x),
\]
with equality in both steps if and only if $C=C_c$.
For $x=0$, only $C=\Pi$ satisfies $\xi(C)=0$, and for $x=1$, the bound $\eta\le1$ is trivial.
Hence
\[
  \eta_+(x)=\max\bigl\{\eta(C):C\in\CC,\ \xi(C)\le x\bigr\},\qquad x\in[0,1],
\]
and likewise $\eta_-(x)=\min\{\eta(C):C\in\CC,\ \xi(C)\le x\}$.

\emph{Step 5: properties of $\eta_\pm$.}
Since $\xi$ is convex and $\eta$ is affine, the representation of $\eta_+$ in Step~4 as a maximum over the convex set $\{\xi\le x\}$ shows that $\eta_+$ is concave on $[0,1]$.
In particular, it is continuous on $(0,1)$.
Continuity at $0$ follows from $0\le\eta_+(x)\le\sqrt{11x/5}$ by Proposition~\ref{prop:sqrt}, and continuity at $1$ from $\eta(C_c)\to1$ as $c\downarrow0$.
Moreover, $\eta_+(x)\ge\frac c2x>0$ for $x=\xi(C_c)\in(0,1)$ by Step~3, and $\eta_+(x)<1$ for $x<1$ because $M$ is the only copula with $\eta=1$.
Together with Step~3, this shows that $\eta_+$ is strictly increasing on $[0,1]$.

\emph{Step 6: filling the fibers.}
Let $x\in(0,1]$ and $e\in[\eta_-(x),\eta_+(x)]$.
By Step~4, there are copulas $C_\pm$ with $\xi(C_\pm)\le x$ and $\eta(C_\pm)=\eta_\pm(x)$.
A suitable mixture $\widetilde C_0$ of $C_-$ and $C_+$ satisfies $\eta(\widetilde C_0)=e$ and, by convexity of $\xi$, also $\xi(\widetilde C_0)\le x$.
If $\xi(\widetilde C_0)<x$, we increase $\xi$ without changing $\eta$ by mixing with a completely dependent copula.
By Lemma~\ref{lem:shuffle-values}, there is $s\in[0,1]$ with $\eta(S_s)=1-2s^4=e$, and $\xi(S_s)=1$ because $S_s$ is supported on the graph of a function, see \cite[Thm.~1.1]{chatterjee2020}.
The mixtures $\widetilde C_t\coloneqq(1-t)\widetilde C_0+tS_s$ satisfy $\eta(\widetilde C_t)=e$ for all $t\in[0,1]$.
By~\eqref{eq:xi-def}, the function $t\mapsto\xi(\widetilde C_t)$ is a quadratic polynomial with $\xi(\widetilde C_0)\le x\le\xi(\widetilde C_1)$, so the intermediate value theorem yields $t$ with $\xi(\widetilde C_t)=x$.
Hence $(x,e)$ belongs to the region.
Together with Step~4 and the fiber $\{0\}$ at $x=0$, this proves the claimed representation of the region.
\qed

\subsection{Comparison with Spearman's rho and Blest's coefficient}

We close the section by relating Theorem~\ref{thm:xi-eta} and Proposition~\ref{prop:sqrt} to the corresponding results for Spearman's rho and for Blest's coefficient itself.

\begin{remark}
\label{rem:comparison}
\leavevmode
\begin{enumerate}[label=(\alph*),leftmargin=2em]
\item
The proof of Proposition~\ref{prop:sqrt} applies verbatim to every coefficient of the form
\[
  \theta(C)=k_0-6\int_0^1\!\int_0^1a_\theta(x,z)\,\kappa(x,z)\de x\de z
\]
with a constant $k_0$ and a bounded weight $a_\theta$ that is Lipschitz in $z$, such that $\theta(\Pi)=0$.
It yields $|\theta(C)|\le\sqrt{\sigma^2\xi(C)}$ with $\sigma^2\coloneqq6\int_0^1\mathrm{Var}(a_\theta(X,z))\de z$ for uniform $X$.
The computation leading to~\eqref{eq:kappa-rep} shows that $\rho$ and $\nu$ are of this form with the weights $2x$ and $2x^2$, respectively.
Since $\mathrm{Var}(2X)=\frac13$ and $\mathrm{Var}(2X^2)=\frac{16}{45}$, we obtain
\[
  |\rho(C)|\le\sqrt{2\,\xi(C)}
  \qquad\text{and}\qquad
  |\nu(C)|\le\sqrt{\tfrac{32}{15}\,\xi(C)}.
\]
The construction in the proof of Proposition~\ref{prop:sqrt} shows that both constants are sharp as $\xi\to0$, in accordance with the explicit boundary curves in \cite{ansari2026exact} and \cite{rockel2026exact}.
\item
The exact $(\xi,\nu)$-region was determined in \cite{rockel2026exact}, and it differs from the $(\xi,\eta)$-region in two respects.
First, it is symmetric about $\nu=0$.
Indeed, the reflection $C\mapsto C^\perp$ changes the sign of $\nu$ by \cite[Lem.~2.2]{rockel2026blest}, and it leaves $\xi$ unchanged because $\partial_1C^\perp(u,v)=1-\partial_1C(u,1-v)$ and the expansion in the proof of Lemma~\ref{lem:kappa} applies.
The fibers of the $(\xi,\eta)$-region, in contrast, are not symmetric about zero by~\eqref{eq:eta-expansion}.
Second, comparing the sharp constants in Proposition~\ref{prop:sqrt} and in~(a) shows that $\eta_+(x)>\max\{\nu(C):\xi(C)\le x\}$ for all sufficiently small $x>0$.
The $(\xi,\eta)$-region is therefore not contained in the $(\xi,\nu)$-region, in contrast to the inclusion of the $(\rho,\eta)$-region in the $(\rho,\nu)$-region from Remark~\ref{rem:rho-nu}.
\item
The last observation is possible because $\xi$ is not symmetric.
For small $x$, the transpose of $C_c^+$ satisfies $\nu((C_c^+)^\top)=2\eta_+(x)-\nu(C_c^+)>\max\{\nu(C):\xi(C)\le x\}$, and hence $\xi((C_c^+)^\top)>x$.
The extremal copulas of the two regions also differ structurally: the conditional distribution functions of the extremal copulas for $\nu$ in \cite{rockel2026exact} are clamped functions of $(1-u)^2$, whereas those of $C_c^\pm$ are clamped functions of the rank sum $u+v$.
\item
Numerically, the two regions are nevertheless close.
The upper boundary of the $(\xi,\eta)$-region lies above that of the $(\xi,\nu)$-region over the whole range, but by less than $2\cdot10^{-3}$.
The lower boundary of the $(\xi,\eta)$-region lies below the lower boundary of the $(\xi,\nu)$-region near $\xi=0$ and above it by up to about $0.018$ for moderate values of $\xi$.
\item
Since $\eta(C^\top)=\eta(C)$, Theorem~\ref{thm:xi-eta} and Proposition~\ref{prop:sqrt} remain valid with $\xi(C^\top)$ in place of $\xi(C)$, with the extremal copulas $(C_c^\pm)^\top$.
\end{enumerate}
\end{remark}

\section{Tail dependence and open problems}
\label{sec:discussion}

\subsection{Tail dependence}

The lower and upper tail dependence coefficients of a copula $C$ are
\[
  \lambda_L(C)\coloneqq\lim_{t\downarrow0}\frac{C(t,t)}{t},
  \qquad
  \lambda_U(C)\coloneqq\lim_{t\uparrow1}\frac{1-2t+C(t,t)}{1-t},
\]
see \cite[Sec.~5.4]{nelsen2006introduction}, defined for those copulas for which the limits exist.
They describe the dependence in the extreme corners only, whereas $\eta$ integrates the copula against a weight that is bounded away from zero on the lower left part of the square.
Accordingly, the tail dependence coefficients carry no information about $\eta$ apart from the two extreme cases $M$ and $W$.

\begin{proposition}[Tail dependence]
\label{prop:tail}
Let $\lambda\in\{\lambda_L,\lambda_U\}$ and let $\CC_\lambda$ be the set of copulas for which $\lambda(C)$ exists.
Then
\[
  \bigl\{(\lambda(C),\eta(C)):C\in\CC_\lambda\bigr\}=\bigl([0,1]\times(-1,1)\bigr)\cup\bigl\{(0,-1),(1,1)\bigr\}.
\]
\end{proposition}

\begin{proof}
\emph{The extreme values of $\eta$.}
For every copula $C$, the functions $M-C$ and $C-W$ are nonnegative and continuous, and the weight $2-u-v$ in~\eqref{eq:eta-def} is positive on $[0,1)^2$.
Hence $\eta(C)=1$ holds only for $C=M$, and $\eta(C)=-1$ only for $C=W$.
Since $\lambda_L(M)=\lambda_U(M)=1$ and $\lambda_L(W)=\lambda_U(W)=0$, the region contains no further points with $|e|=1$.

\emph{Interior points.}
Let $l\in[0,1]$ and $e\in(-1,1)$, and consider $\lambda=\lambda_L$.
For $\varepsilon\in(0,1)$ and $\widetilde C\in\CC$, let $O_{\varepsilon,\widetilde C}$ be the ordinal sum, see \cite[Sec.~3.2.2]{nelsen2006introduction}, of the copula $lM+(1-l)\Pi$ on $[0,\varepsilon]^2$ and of $\widetilde C$ on $[\varepsilon,1]^2$.
For $t\le\varepsilon$, we have
\[
  O_{\varepsilon,\widetilde C}(t,t)=\varepsilon\Bigl(l\,\frac t\varepsilon+(1-l)\frac{t^2}{\varepsilon^2}\Bigr),
\]
and hence $\lambda_L(O_{\varepsilon,\widetilde C})=l$, whatever $\widetilde C$ is.
The ordinal sum is affine in $\widetilde C$, hence so is $\eta(O_{\varepsilon,\widetilde C})$.
Moreover, $O_{\varepsilon,W}\to W$ and $O_{\varepsilon,M}\to M$ uniformly as $\varepsilon\downarrow0$.
Since $\eta$ is continuous with respect to uniform convergence, we can fix $\varepsilon$ with $\eta(O_{\varepsilon,W})<e<\eta(O_{\varepsilon,M})$.
The intermediate value theorem along the segment $\widetilde C=(1-t)W+tM$, $t\in[0,1]$, then yields a copula with $\lambda_L=l$ and $\eta=e$.
For $\lambda_U$, the block $lM+(1-l)\Pi$ is placed on $[1-\varepsilon,1]^2$ instead.
\end{proof}

The same construction with blocks at both corners shows that the joint range of $(\lambda_L,\lambda_U,\eta)$ is $([0,1]^2\times(-1,1))\cup\{(0,0,-1),(1,1,1)\}$.

\subsection{Open problems}

Together with the $(\nu,\eta)$-region of \cite{rockel2026blest}, Theorems~\ref{thm:rho-eta}, \ref{thm:beta-eta} and~\ref{thm:xi-eta} and Proposition~\ref{prop:tail} determine the exact regions of $\eta$ with Spearman's rho, Blomqvist's beta, Chatterjee's rank correlation, Blest's coefficient and the tail dependence coefficients.
The regions of $\eta$ with Spearman's footrule, Gini's gamma, Kendall's tau, the van der Waerden coefficient and the measures of Schweizer--Wolff and Hoeffding remain open.
We comment on the first three.

Spearman's footrule, Gini's gamma and the van der Waerden coefficient are affine in the copula, just like Spearman's rho.
The corresponding regions are therefore convex, and their supporting lines are again optimal transport problems.
The costs, however, are kinked along one or both diagonals of the square, or they involve the normal quantile function.
Linear programs over checkerboard copulas indicate that parts of these boundaries are attained by ordinal sums of increasingly many cyclic shifts, or by couplings supported on curved graphs, so we do not expect boundaries as simple as those of Theorem~\ref{thm:rho-eta}.

For Kendall's tau, the shuffles $S_s$ are not extremal.
To see this, consider ordinal sums of copies of $M$ and $W$.
If the countermonotone blocks have lengths $\ell_i$, then two points fall into the same countermonotone block with probability $\sum_i\ell_i^2$, and all other pairs are concordant, so that $\tau=1-2\sum_i\ell_i^2$.
A countermonotone block on $[t,t+\ell]^2$ lowers $\eta$ by $2\ell^3(2-2t-\ell)$ compared with $M$.
Now let $s\in(0,1)$ and $0<\ell_1<\min\{1-s,s/\sqrt2\}$, and replace the block $[1-s,1]^2$ of $S_s$ by two adjacent countermonotone blocks of lengths $\ell_1$ and $\ell_2\coloneqq\sqrt{s^2-\ell_1^2}>\ell_1$, the longer one at the end.
This preserves $\tau$, and it raises $\eta$ by
\[
  2s^4-2\ell_2^4-2\ell_1^3(\ell_1+2\ell_2)=4\ell_1^2\ell_2(\ell_2-\ell_1)>0.
\]
As for the $(\tau,\rho)$-region of \cite{schreyer2017exact}, we expect the upper boundary of the $(\tau,\eta)$-region to be attained by ordinal sums with several countermonotone blocks.

\appendix

\section{Proof of Proposition~\ref{prop:closed}}
\label{app:closed}

We work in the reflected coordinates of Section~\ref{sec:xi} and write $y=x+z$ for the rank sum at the threshold $z$.
The proof consists in solving the normalization $\int_0^1\kappa_c^\pm(x,z)\de x=z$ in each configuration of the active interval, in computing the integrals in~\eqref{eq:kappa-rep} level by level, and in integrating over the thresholds.

For the upper family $\kappa_c^+$, we put $B=B_c(z)$ and $A=\sqrt{\max\{B^2-c,0\}}$ as in the proof of Lemma~\ref{lem:families}.
The kernel equals $1$ for $y\le A$, it equals $(B^2-y^2)/c$ for $A<y<B$, and it vanishes for $y\ge B$.
Four configurations can occur, according to whether the active interval $(A,B)$ extends beyond the left endpoint $z$ or the right endpoint $1+z$ of the range of the rank sum.
In each of them, the normalization can be solved by a suitable parametrization.
\begin{enumerate}[label=(U\arabic*),leftmargin=2.8em]
\item \emph{$A<z$ and $B\le1+z$.}
The kernel is active at $x=0$ and vanishes before $x=1$, and the normalization reads
\[
  (B-z)^2(2B+z)=3cz.
\]
With the parameter $s\coloneqq B-z$, the $x$-length of the support, it is solved by $z=\frac{2s^3}{3(c-s^2)}$.
\item \emph{$A\ge z$ and $B\le1+z$.}
The kernel equals $1$ near $x=0$ and vanishes before $x=1$.
The normalization reads
\[
  2z=A+\frac{(2B+A)(B-A)}{3(A+B)},
\]
and with $\vartheta\coloneqq(A+B)/\sqrt c$, it is solved by
\[
  A=\tfrac{\sqrt c}{2}(\vartheta-\vartheta^{-1}),
  \qquad
  B=\tfrac{\sqrt c}{2}(\vartheta+\vartheta^{-1}),
  \qquad
  z=\sqrt c\,\Bigl(\frac\vartheta4+\frac{1}{12\vartheta^3}\Bigr).
\]
\item \emph{$A\ge z$ and $B>1+z$.}
The kernel equals $1$ near $x=0$ and is active at $x=1$.
The normalization reads
\[
  (1+z-A)^2(1+z+2A)=3c(1-z),
\]
and with $r\coloneqq1+z-A$, it is solved by $z=\frac{3c-3r^2+2r^3}{3(c+r^2)}$.
\item \emph{$A<z$ and $B>1+z$.}
The kernel is active on all of $[0,1]$, and the normalization reads
\[
  B^2=z^2+(1+c)z+\tfrac13.
\]
\end{enumerate}
As $z$ increases from $0$ to $1$, the kernel passes through three of these configurations.
For $c\ge c_+$, the order is (U1) for $s\in(0,1]$, then (U4) for $z\in[\frac{2}{3(c-1)},\frac{3c-1}{3(c+1)}]$, then (U3) for $r\in(0,1]$.
For $c\le c_+$, which is the range of~\eqref{eq:up-small}, the order is (U1) for $s\in(0,\sqrt{(2\sqrt3-3)c}]$, then (U2) for $\vartheta\in[\sqrt{c_+},P]$ with $P$ from~\eqref{eq:P-def}, then (U3) for $r\in(0,\sqrt c/P]$.
At the branch point $c=c_+$, the middle configuration degenerates to a single threshold, which explains the change of form of the closed expressions.
The map $P\mapsto c$ defined by~\eqref{eq:P-def} is a decreasing bijection from $[\sqrt{c_+},\infty)$ onto $(0,c_+]$.

For the lower family $\kappa_c^-$, we put $L=L_c(z)$, $D=\sqrt{\max\{L,0\}}$ and $E=\sqrt{L+c}$.
The kernel vanishes for $y\le D$, it equals $(y^2-L)/c$ for $D<y<E$, and it equals $1$ for $y\ge E$.
\begin{enumerate}[label=(L\arabic*),leftmargin=2.8em]
\item \emph{$D\ge z$ and $E>1+z$.}
The normalization reads
\[
  (1+z-D)^2(1+z+2D)=3cz,
\]
and with $t\coloneqq1+z-D$, it is solved by $z=\frac{t^2(3-2t)}{3(c-t^2)}$.
\item \emph{$D\ge z$ and $E\le1+z$.}
The normalization reads
\[
  2(D^2+DE+E^2)=3(D+E),
\]
which does not involve $z$.
Together with $E^2-D^2=c$, it determines the constants
\[
  D=\tfrac{\sqrt c}{2}(p-p^{-1}),
  \qquad
  E=\tfrac{\sqrt c}{2}(p+p^{-1}),
\]
where $p$ is given by~\eqref{eq:p-def}.
This configuration covers the thresholds $z\in[E-1,D]$.
\item \emph{$D<z$ and $E\le1+z$.}
The normalization reads
\[
  (E-z)^2(2E+z)=3c(1-z),
\]
and with $\mu\coloneqq E-z$, it is solved by $z=\frac{3c-2\mu^3}{3(c+\mu^2)}$.
\item \emph{$D<z$ and $E>1+z$.}
The kernel is active on all of $[0,1]$, and the normalization reads
\[
  L=z^2+(1-c)z+\tfrac13.
\]
\end{enumerate}
For $c\ge c_-$, the order is (L1) for $t\in(0,1]$, then (L4) for $z\in[\frac{1}{3(c-1)},\frac{3c-2}{3(c+1)}]$, then (L3) for $\mu\in(0,1]$.
For $c\le c_-$, the order is (L1) for $t\in(0,\sqrt c/p]$, then (L2), then (L3) for $\mu\in(0,\sqrt c/p]$.
The map $p\mapsto c$ defined by~\eqref{eq:p-def} is a decreasing bijection from $[\sqrt{c_-},\infty)$ onto $(0,c_-]$.
In configuration (L2), the kernel $\kappa_c^-(x,z)$ depends on $x$ and $z$ only through the rank sum $x+z$, so the density of $C_c^-$ is a function of $u+v$ on the corresponding band.

It remains to evaluate the integrals.
In each configuration, the level integrals
\[
  \int_0^1\kappa_c^\pm(x,z)^2\de x
  \qquad\text{and}\qquad
  \int_0^1a(x,z)\,\kappa_c^\pm(x,z)\de x
\]
are integrals of polynomials in $y$ over intervals whose endpoints are among $z$, $1+z$, $A$, $B$, $D$ and $E$.
They are therefore polynomials in $z$ and the configuration parameter.
Changing the variable of integration from $z$ to that parameter, the integrals over $z$ in~\eqref{eq:kappa-rep} become integrals of rational functions of the parameter, with denominators that are powers of $c-s^2$, $c+r^2$, $c-t^2$ or $c+\mu^2$, or monomials in $\vartheta$.
Their antiderivatives consist of rational functions, logarithms and arctangents.
Summing the contributions of the three configurations yields~\eqref{eq:up-large}--\eqref{eq:lo-small}.
We carried out these integrations with computer algebra and checked the results against a direct numerical evaluation of~\eqref{eq:kappa-rep} for the numerically computed kernels.
The two branches of each curve agree at $c=c_\pm$, and the expansions in Proposition~\ref{prop:closed} follow from~\eqref{eq:up-large} and~\eqref{eq:lo-large} by expanding in powers of $c^{-1/2}$.

\bibliographystyle{plainnat}
\bibliography{references}

\end{document}